\documentclass{article}
\usepackage{amsrefs}
\usepackage{amssymb}
\usepackage{amsmath}
\usepackage{amsthm}
\usepackage{graphicx}
\DeclareMathOperator{\Supp}{Supp}\DeclareMathOperator{\Div}{Div}
\DeclareMathOperator{\Prin}{Prin}\DeclareMathOperator{\Spec}{Spec}
\DeclareMathOperator{\rk}{rk}\DeclareMathOperator{\val}{val}
\numberwithin{equation}{section}
\newtheorem{theorem}{Theorem}[section]
\newtheorem{definition}[theorem]{Definition}
\newtheorem{example}[theorem]{Example}
\newtheorem{proposition}[theorem]{Proposition}
\newtheorem{lemma}[theorem]{Lemma}
\newtheorem{definition/construction}[theorem]{Definition/Construction}
\newtheorem{construction}[theorem]{Construction}
\theoremstyle{remark}
\newtheorem*{remark}{Remark}

\author{Marc Coppens\footnote{Katholieke Hogeschool Kempen, Departement IBW,
Kleinhoefstraat 4, B-2440 Geel, Belgium; K.U.Leuven, Departement of Mathematics,
Celestijnenlaan 200B, B-3001 Leuven, Belgium; email: marc.coppens@khk.be}}
\title{Linear systems on graphs with a real structure}
\date{}

\begin{document}
\maketitle \noindent

\begin{abstract}

A degeneration of a smooth projective curve to a strongly stable
curve gives rise to a specialization map from divisors on curves to
divisors on graphs. In this paper we show that this specialization
behaves well under the presence of real structures. In particular we
study real linear systems on graphs with a real structure and we
prove results on them comparable to results in the classical theory
of real curves. We also consider generalizations to metric graphs
and tropical curves

\end{abstract}

\section{Introduction}\label{section1}

In \cite{ref1} one introduces a theory of linear systems on graphs
which is very similar to the theory of linear systems on smooth
complex projective curves. In \cite{ref2} one shows a relation by
means of specialization of a smooth projective curve to a strongly
semistable curve. Associated to this strongly semistable curve there
is a dual graph $G$ and divisors on the geometric generic smooth
fiber give rise to divisors on the metric graph $T$ associated to
$G$. This specialization preserves linear equivalence on curves and
on graphs. It is the aim of this paper to show that this
specialization is compatible with real structures. This means in
case there is a real structure on the family of curves there is a
real structure on the graph and on the metric graph and real
divisors on curves specialize to real divisors on the graph.

In Section~\ref{section2} we introduce the notion of a graph $G$
with a real structure (Definition \ref{definitie1}). We define an
associated real locus graph $G(\mathbb{R})$ (Definition
\ref{definition5}). This is similar to the notion of the real locus
$X(\mathbb{R})$ of a real curve $X$. In this classical case there
are strong constraints on the number of connected components $s(X)$
of $X(\mathbb{R})$, the genus $g(X)$ and a number $a(X)$ determined
by the number of connected components of $X(\mathbb{C})\setminus
X(\mathbb{R})$. A graph $G$ has a genus $g(G)$ and we define similar
numbers $s(G)$ (Formula \ref{formule1}) and $a(G)$ (Definition
\ref{definition6}). The definition of $s(G)$ not only takes into
account the number of connected components of $G(\mathbb{R})$ but
also their genus (in the curve case all components of
$X(\mathbb{R})$ are circles). In the graph case we prove there are
similar but stronger constraints (Theorem \ref{theorem1}). We prove
that all cases satisfying those constraints do occur (Example
\ref{voorbeeld1}).

As mentioned before, there is a theory of divisors and linear
systems on graphs similar to the theory on complex curves. In case
the graph has a real structure one can consider real divisors and
real linear systems. In Section~\ref{section3} we prove results
similar to basic properties of real linear systems on real curves.
In particular we prove that the parity of the degree of the
restriction of a real divisor $D$ to a connected component of
$G(\mathbb{R})$ does not change under linear equivalence (Theorem
\ref{theorem2}). We introduce an M-graph as being a graph satisfying
$s(G)=g(G)+1$ (Definition \ref{definitie8}; this is similar to the
concept of an M-curve). We prove that each real effective divisor on
an M-graph is linearly equivalent to a totally real effective
divisor (Theorem \ref{theorem3}). A similar strong fact does not
hold in the case of M-curves. In case $G(\mathbb{R})$ has $g(G)+1$
connected components we say that $G$ is a strong M-graph (Definition
\ref{definitie10}) and we prove $G$ has a real linear system $g^1_2$
in that case (Proposition \ref{propositie3}).

In order to study the compatibility of the above mentioned
specialization with real structures we need to consider metric
graphs with a real structure. This is done in
Section~\ref{section4}. We generalize the results on real linear
systems on graphs to the context of metric graphs with a real
structure. We also make further generalizations to tropical curves.

Finally in Section~\ref{section5} we study the behavior of the real
structure under specialization. Although at the end this follows the
arguments of \cite{ref2} quite closely (see the proof of Theorem
\ref{theorem6}), we need to take care about this behavior of the
real structure in those arguments. In particular for the generic
fiber of the degeneration the base field is a real field and we use
the real closure instead of the the algebraic closure to define real
divisors on the generic fiber. Therefore we need to pay some
attention to extensions of real fields in those arguments (Lemmas
\ref{lemma3}, \ref{lemma4} and \ref{lemma5}). Also we need to take
care about the existence of real structures on certain
desingularizations occurring in those arguments (Lemma \ref{lemma6})
and we need to check that this real structure is compatible to the
real structure on the metric graph $T$ associated to $G$
(Proposition \ref{proposition8}). Finally we prove that all graphs
with a real structure can be obtained from such a degeneration
(Proposition \ref{proposition9}).

\section{Graphs with a real structure}\label{section2}

We recall some generalities on finite graphs (as a general
reference, see e.g. \cite{ref10}). A \emph{finite graph} $G$ is
defined by a finite set of vertices $V(G)$ and a finite set of edges
$E(G)$ together with an incidence function that associates to an
edge $e$ a subset $\psi (e)\subset V(G)$ containing one or two
vertices. In case it contains only one vertex $v$ we say $e$ is a
\emph{loop} at $v$. When we write $\psi (e)=\{ v,w\}$ it is possible
that $v=w$ (i.e. $e$ could be a loop). Vertices of $\psi (e)$ are
called the \emph{ends} of $e$. We write $v(G)$ (resp. $e(G)$) to
denote the number of vertices (resp. edges) of $G$.

A \emph{walk} $\Gamma$ in $G$ is a sequence $v_0 e_0 v_1 e_1 v_2
\cdots v_{n-1} e_{n-1} v_n$ such that $e_i\in E(G)$ with $\psi
(e_i)=\{ v_i, v_{i+1}\}$ for $0\leq i\leq n-1$. We call $v_0$ and
$v_n$ the \emph{ends} of the walk $\Gamma$ and we say $\Gamma$
connects $v_0$ to $v_n$. We call $e_i$ ($0\leq i\leq n-1$) the edges
of $\Gamma$ and $v_i$ ($1\leq i\leq n-1$) the internal vertices of
$\Gamma$. A walk $\Gamma$ is called closed if the ends $v_0$ and
$v_n$ of $\Gamma$ are equal. A graph $G$ is \emph{connected} if for
each two different vertices $v_1, v_2$ of $G$ there exists a walk
having ends $v_1$ and $v_2$.

A \emph{subgraph} $G'$ of $G$ is defined by subsets $V(G')\subset
V(G)$ and $E(G')\subset E(G)$ such that for $e\in E(G')$ one has
$\psi (e)\subset V(G')$. If $G'$ and $G''$ are two subgraphs of $G$
then $G'\cup G''$ is the subgraph defined by $V(G')\cup
V(G'')\subset V(G)$ and $E(G')\cup E(G'')\subset E(G)$. A subgraph
$C$ of $G$ is called a \emph{cycle} if all vertices and edges of $C$
are the vertices and edges of a closed walk such that the end
vertices are the only vertices being equal in that walk. A
\emph{connected component} of a graph $G$ is a subgraph $G'$ such
that, each subgraph $G''\neq G'$ containing $G'$ is not connected. A
graph $G$ has finitely many connected components $G_1; \cdots; G_c$.
We write $c=c(G)$ to denote the number of connected components of
$G$. The \emph{genus} of the graph $G$ is defined by

\[
g(G)=c(G)+e(G)-v(G).
\]

From now on in this paper $G$ denotes a connected graph.

\begin{definition}\label{definitie1}
Let $G$ be a finite graph. A \emph{real structure} on $G$ is defined
by an involution $\iota :G\rightarrow G$ (hence $\iota \circ \iota$
is the identity on $G$). More concretely, it is defined by
involutions $\iota _V : V(G)\rightarrow V(G)$ and $\iota _E :
E(G)\rightarrow E(G)$ such that for each $e\in E(G)$ one has $\psi
(\iota _E (e))=\iota _V (\psi (e))$.
\end{definition}

\noindent From now on we assume $G$ is a graph with a fixed real
structure. For $v\in V(G)$ (resp. $e\in E(G)$) we write
$\overline{v}$ (resp. $\overline{e}$) instead of $\iota _V(v)$
(resp. $\iota _E(e)$).

\begin{definition}\label{definitie2}
A vertex $v\in V(G)$ is called a \emph{real vertex} if
$v=\overline{v}$. In case $v$ is not a real vertex $\overline{v}$ is
called the \emph{conjugated vertex}. In that case the symbol
$v+\overline{v}$ (a real divisor on $G$ according to the terminology
of Section \ref{section3}) is called a \emph{non-real vertex pair}
of $G$.
\end{definition}

\noindent We write $V_{\mathbb{R}}(G)$ to denote the subset of
$V(G)$ of real vertices of $G$.

\begin{definition}\label{definitie3}
An edge $e\in E(G)$ is called a \emph{real edge} if
$e=\overline{e}$. A real edge $e$ is called an \emph{isolated real
edge} if $\psi (e)\nsubseteq V_{\mathbb{R}}(G)$, otherwise we say
$e$ is a \emph{non-isolated real edge}. In case $e$ is not a real
edge $\overline{e}$ is called the \emph{conjugated edge} of $e$.
\end{definition}

We write $E_{\mathbb{R}}(G)$ (resp. $E^0_{\mathbb{R}}(G)$) to denote
the subset of $E(G)$ of real edges (resp. non-isolated real edges)
of $G$.

\begin{definition}\label{definitie4}
A subgraph $G'$ of $G$ is called a \emph{real subgraph} if the
subsets $V(G')$ of $V(G)$ and $E(G')$ of $E(G)$ are invariant under
$\iota_V$ and $\iota_E$.
\end{definition}

\noindent If $G'$ is a real subgraph of $G$ then the restrictions of
$\iota_V$ to $V(G')$ and $\iota_E$ to $E(G')$ induce a real
structure on $G'$. We always consider this real structure on $G'$,
hence a real subgraph of $G$ is always considered as a graph with a
real structure.

\begin{definition}\label{definition5}
The subgraph of $G$ having $V_{\mathbb{R}}(G)$ (resp.
$E^0_{\mathbb{R}}(G)$) as its vertex set (resp. edge set) is called
the \emph{real locus graph} of $G$. It is denoted by
$G(\mathbb{R})$.
\end{definition}

\noindent Clearly $G(\mathbb{R})$ is a real subgraph of $G$ and the
identity defines the real structure on $G(\mathbb{R})$. The number
of connected components of $G(\mathbb{R})$ is denoted by $s'(G)$. We
denote the components by $G(\mathbb{R})_i$ ($1\leq i\leq s'(G)$).
Let $e^i(G)$ be the number of isolated real edges of $G$. Now we
introduce the number $s(G)$ that corresponds to the number of
connected components of $X(\mathbb{R})$ in case $X$ is a real curve
as follows:

\begin{equation}\label{formule1}
    s(G)=e^i(G)+\sum_{i=1}^{s'(G)} \bigl( g(G(\mathbb{R})_i)+1
    \bigr).
\end{equation}

\begin{remark}
The definition of this number $s(G)$ is in accordance with the
degeneration studied in Section \ref{section5}. Consider the graph
$G$ as the dual graph of a total degeneration of a smooth real curve
of genus $g(G)$ (see Proposition \ref{proposition9}). In such
degeneration one should compare $G(\mathbb{R})_i$ to a real
algebraic curve of genus $g(G(\mathbb{R})_i)$, the real locus of
such curve has at most $g(G(\mathbb{R})_i)+1$ connected components.
On the other hand, an isolated real edge corresponds to an isolated
real singular point of a real curve. Such isolated real singular
point can deform to a connected component of the real locus of a
smooth curve. In this way $s(G)$ corresponds to the maximal number
of connected components of the real locus of a smooth real curve
having a total degeneration with dual graph $G$ with its real
structure.
\end{remark}

Next we define the number $a(G)$ similar to the number $a(X)$ of a
real curve $X$ defined by $a(X)=1$ if $X(\mathbb{C})\setminus
X(\mathbb{R})$ is connected and $a(X)=0$ otherwise.

\begin{definition}\label{definition6}
Let $G$ be a graph with a real structure. In case there exists a
non-real vertex $v$ and a walk $\Gamma$ in $G$ with ends $v$ and
$\overline{v}$ such that $\Gamma$ contains neither a real vertex,
nor a real edge, we put $a(G)=1$. Otherwise we put $a(G)=0$.
\end{definition}

The following theorem gives conditions between the integers $g(G)$,
$a(G)$ and $s(G)$.

\begin{theorem}\label{theorem1}
Let $G$ be a graph with a real structure. One has $s(G) \equiv
g(G)+1 \pmod{2}$ and $0\leq s(G)\leq g(G)+1$. Moreover $s(G)\neq
g(G)+1$ in case $a(G)=1$ and $s(G)\neq 0$ in case $a(G)=0$.
\end{theorem}

\begin{remark}
In the case of a real smooth curve $X$ the condition $s(X) \equiv
g(X)+1 \pmod{2}$ need not hold if $a(X)=1$. Of course, in the case
of graphs, this condition implies $s(G)\leq g(G)-1$ in case
$a(G)=1$.
\end{remark}

Before giving the proof of Theorem \ref{theorem1} we introduce some
notations we use in that proof. Let $G$ be a graph and let $W$ be a
subset of $V(G)$. Then $G[W]$ is the subgraph \emph{induced by} $W$.
This means, the set of vertices of $G[W]$ is equal to $W$ and the
edges of $G[W]$ are the edges $e$ of $G$ such that $\psi (e)\subset
W$. Let $S$ be a subset of $E(G)$ then $G[S]$ is the subgraph of $G$
whose set of edges is equal to $S$ and the set of vertices of $G[S]$
is the union $\bigcup_{e\in S}\psi (e)$. The graph $G\setminus S$ is
the subgraph of $G$ obtained by deleting the edges belonging to $S$
(hence its set of vertices is equal to $V(G)$). Finally $G(S)$ is
the graph obtained from $G$ by contracting each edge not contained
in $S$. Since $G$ is connected it is obtained from $G[S]$ by
identifying two vertices $v$ and $v'$ if and only if there is a walk
$P$ in $G\setminus S$ having ends $v$ and $v'$. For a subset $S$ of
$E(G)$ we use the equation

\begin{equation}\label{formule2}
    g(G)=g(G(S))+g(G\setminus S).
\end{equation}

\noindent Clearly, even if $G$ is a connected graph, the graph
$G\setminus S$ need not be connected. In particular it can have
isolated vertices as some connected components. Such components do
not contribute to $g(G\setminus S)$.

\begin{proof}[Proof of Theorem \ref{theorem1}]\hfill

\begin{proof}[Step 1]
Reduction to the case that $G$ has no isolated real edges.
\renewcommand{\qedsymbol}{ }
\end{proof}

\noindent Assume $e$ is an isolated real edge of $G$ with $\psi
(e)=\{ v;\overline{v}\}$ for some non-real vertex pair
$v+\overline{v}$. Let $G'$ be the graph with real structure obtained
from $G$ by first deleting the edge $e$ (this is a real subgraph of
$G$) and then adding a new vertex $v_e$ and two new edges $e'$ and
$e''$ with $\psi (e')=\{ v; v_e\}$ and $\psi (e'')=\{ \overline{v};
v_e\}$ and putting $\iota (v_e)=v_e$ and $\iota (e')=e''$, $\iota
(e'')=e'$. Instead of the isolated real edge $e$ this graph has one
more component $\{ v_e\}$ of $G'(\mathbb{R})$. Clearly $g(G')=g(G)$,
$s(G')=s(G)$ and $a(G')=a(G)$. Continuing in the way this shows it
is enough to prove the case without isolated real edges. So from now
on we assume $e^i(G)=0$.

\begin{figure}[h]
\begin{center}
\includegraphics[height=4 cm]{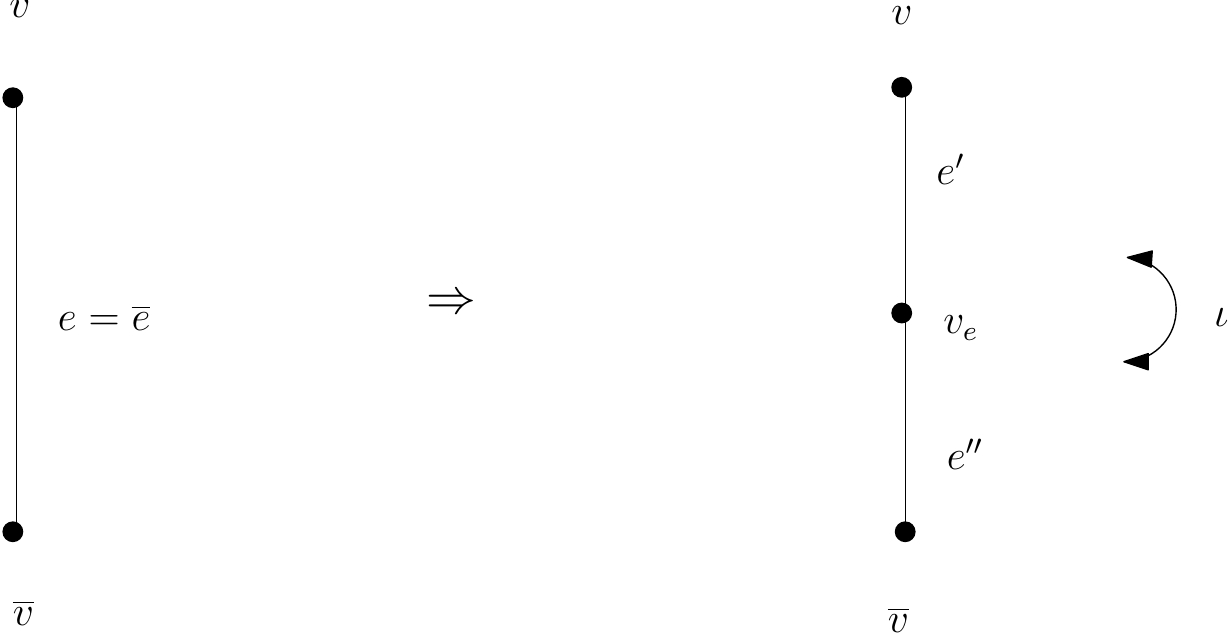}
\caption{Replacing an isolated real edge by a real vertex.}
\end{center}
\end{figure}

\begin{proof}[Step 2]
The case $G(\mathbb{R})$ is empty.
\renewcommand{\qedsymbol}{ }
\end{proof}

\noindent In case $G(\mathbb{R})$ is empty one has $a(G)=1$ and
$s(G)=0$ and we need to prove that $g(G)$ is odd. Choose a vertex
$v$ of $G$ and a walk $\Gamma$ in $G$ that connects $v$ to
$\overline{v}$. Writing $\Gamma =v_0(=v)e_0v_1 \cdots
v_{n-1}e_{n-1}v_n(=\overline{v})$ we can assume that for each $0\leq
i<j\leq n$ we have $v_j\notin \{v_i;\overline{v_i}\}$, otherwise we
replace $\Gamma$ by a shorter walk. It follows that $\Gamma \cup
\overline{\Gamma}$ is a cycle $C$ with $\overline{C}=C$. Let $e$ be
an edge of $C$ with $v\in \psi (e)$ and consider the graph
$G'=G\setminus \{ e; \overline{e}\}$. Clearly $C\setminus \{e;
\overline{e}\}$ has exactly two connected components $C_1$ and $C_2$
with $C_2=\overline{C_1}$. Let $w\in V(G)\setminus V(C)$. Since $G$
is connected there exists $w'\in V(C)$ and a walk $\Gamma_w$ that
connects $w$ to $w'$ such that $\Gamma_w$ contains no edge of $C$.
In case $w'\in V(C_1)$, $w$ belongs to the connected component $G_1$
of $G'$ containing $C_1$. In that case $\overline{\Gamma_w}$
connects $\overline{w}$ to $\overline{w'}\in V(C_2)$, hence
$\overline{w}$ belongs to the connected component $G_2$ of $G'$
containing $C_2$.

\begin{figure}[h]
\begin{center}
\includegraphics[height=5 cm]{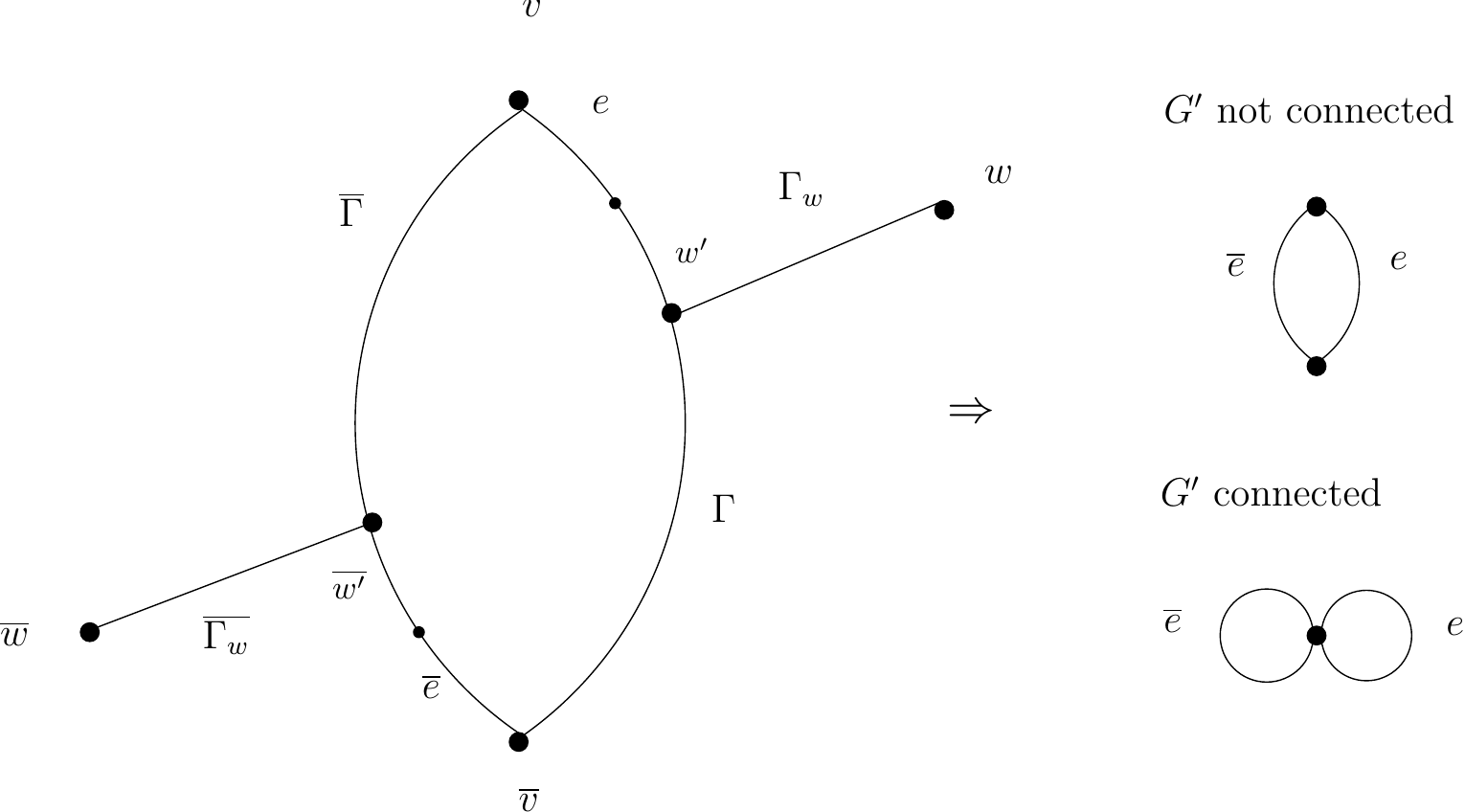}
\caption{Making $G(\{e,\overline{e}\})$.}
\end{center}
\end{figure}

It follows that $G'$ is either connected or it has exactly two
connected components $G_1$ and $G_2$. In case $G'$ is not connected
$G(\{ e;\overline{e}\})$ is a cycle, hence it has genus 1 and
$g(G')=2g(G_1)$. Hence Formula \eqref{formule2} implies
$g(G)=1+2g(G_1)$, in particular $g(G)$ is odd. In case $G'$ is
connected $G(\{ e;\overline{e}\})$ consists of one vertex and two
loops, hence it has genus 2. Hence Formula \eqref{formule2} implies
$g(G) \equiv g(G') \pmod{2}$. Moreover $G'$ is a connected graph
with a  real structure satisfying $s(G')=0$, hence $a(G')=0$. Since
$\sharp (E(G'))<\sharp (E(G))$ we can use induction and assume
$g(G')$ is odd. It follows that $g(G)$ is odd too.

From now on we assume $G(\mathbb{R})$ is not empty.

\begin{proof}[Step 3]
It is enough to prove the theorem in case $a(G)=0$.
\renewcommand{\qedsymbol}{ }
\end{proof}

\noindent Assume the theorem holds in case $a=0$ and assume
$a(G)=1$. There exists a non-real vertex pair $v+\overline{v}$ and a
walk $\Gamma$ that connects $v$ to $\overline{v}$ such that $\Gamma$
contains no real vertex and no real edge. As in the previous step we
can assume $C=\Gamma \cup \overline{\Gamma}$ is a cycle and we take
$e$ and $\overline{e}$ as in the proof of the previous step. Since
$G(\mathbb{R})\neq \emptyset$ there is some real vertex $w$. In case
$G\setminus \{ e,\overline{e}\}$ has two connected components $G_1$
and $G_2$ it would follow that $w\in G_1\cap G_2$ since
$G_2=\overline{G_1}$. This is impossible, hence $G\setminus \{e,
\overline{e}\}$ is a connected graph $G'$ with $G'(\mathbb{R})\neq
\emptyset$. As explained in the previous step this graph satisfies
$s(G)=s(G')$ and $g(G')=g(G)-2$. In case $a(G')=0$, $s(G') \equiv
g(G')+1 \pmod{2}$ and $s(G')\leq g(G')+1$ implies $s(G) \equiv
g(G)+1 \pmod{2}$ and $s(G)\leq g(G)-1$. In case $a(G')=1$ we can
continue this procedure that has to finish because the genus becomes
less.

From now on we can assume $a(G)=0$.

\begin{proof}[Step 4]
The case $G(\mathbb{R})$ is connected.
\renewcommand{\qedsymbol}{ }
\end{proof}

\noindent Consider the induced subgraph $G[V_{\mathbb{R}}(G)]$. Of
course $G(\mathbb{R})$ is a subgraph of $G[V_{\mathbb{R}}(G)]$
having the same set of vertices $V_{\mathbb{R}}(G)$. If $e$ is an
edge in $G[V_{\mathbb{R}}(G)]$ not contained in $G(\mathbb{R})$ then
$\psi (e)\subset V_{\mathbb{R}}(G)$, hence $\psi (e)=\psi
(\overline{e})$ and therefore $\overline{e}$ is another edge in
$G[V_{\mathbb{R}}(G)]$ not contained in $G(\mathbb{R})$. Since
$G(\mathbb{R})$ is connected it follows that
$g(G[V_{\mathbb{R}}(G)])\geq g(G(\mathbb{R}))$ with equality if and
only if $G[V_{\mathbb{R}}(G)]=G(\mathbb{R})$. Moreover we obtain
$g(G[V_{\mathbb{R}}(G)]) \equiv g(G(\mathbb{R})) \pmod{2}$ and since
$g(G(\mathbb{R}))=s(G)-1$ in case $G(\mathbb{R})$ is connected we
obtain

\begin{equation}\label{formule3}
    g(G[V_{\mathbb{R}}(G)]) \equiv s(G)+1 \pmod{2}.
\end{equation}

\noindent This finishes the proof of the theorem in case
$G=G[V_{\mathbb{R}}(G)]$, so we assume that $G\neq
G[V_{\mathbb{R}}(G)]$.

\begin{figure}[h]
\begin{center}
\includegraphics[height=5 cm]{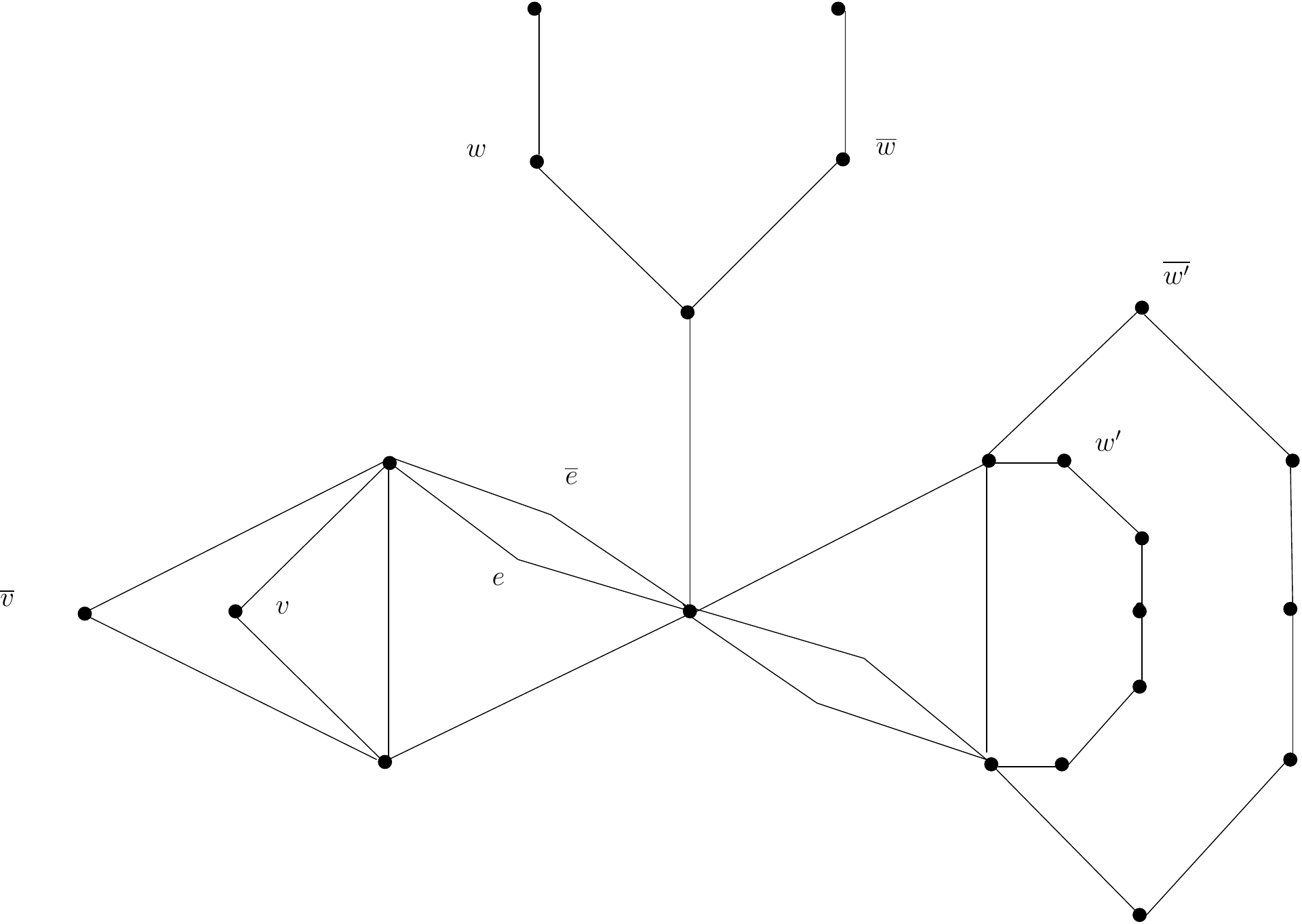}
\caption{Example of a graph $G$ with connected $G(\mathbb{R}).$}
\end{center}
\end{figure}

\begin{figure}[h]
\begin{center}
\includegraphics[height=4 cm]{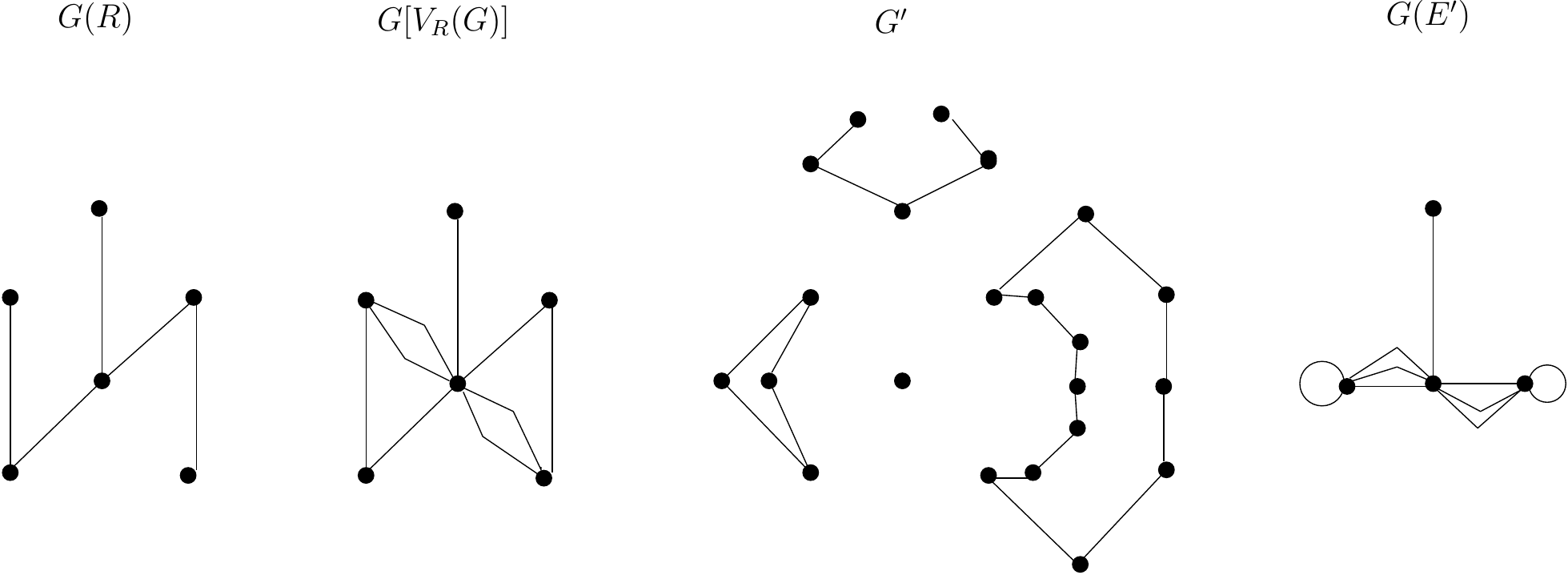}
\caption{Graphs $G'$ and $G(E')$ of the example.}
\end{center}
\end{figure}

Assume $G[V_{\mathbb{R}}(G)]$ is just one real vertex $v$ and
contains no edge. Let $G'$ be the graph obtained from $G$ by adding
a loop $e$ at $v$ and by extending the real structure by putting
$\iota _E(e)=e$. Then $g(G')=g(G)+1$, $a(G')=a(G)$ and
$s(G')=s(G)+1$. Therefore we can assume $G[V_{\mathbb{R}}(G)]$ is
not equal to just one real vertex and no edge. Since
$G[V_{\mathbb{R}}(G)]$ is connected it follows that
$G[V_{\mathbb{R}}(G)]=G[E']$ with $E'$ the set of edges of
$G[V_{\mathbb{R}}(G)]$. Let $G'=G\setminus E'$, it is a real
subgraph of $G$. We are going to use the decomposition (see Formula
\ref{formule2})

\begin{equation}\label{formule4}
    g(G)=g(G(E'))+g(G').
\end{equation}

The graph $G(E')$ is obtained from $G[E']=G[V_{\mathbb{R}}(G)]$ by
identifying vertices $v_1$ and $v_2$ of $V_{\mathbb{R}}(G)$ if and
only if there is a walk $\Gamma$ in $G'$ with ends $v_1$ and $v_2$.
In particular we obtain

\begin{equation}\label{formule5}
    g(G(E'))\geq g(G[V_{\mathbb{R}}(G)])\geq
    g(G(\mathbb{R}))=s(G)-1.
\end{equation}

\noindent Together with Formula \ref{formule4} this already proves
$s(G)\leq g(G)+1$. Moreover we obtain the following claim that will
be used in Section \ref{section3}.

\begin{proof}[Claim]
If $G(\mathbb{R})$ is connected then $s(G)=g(G)+1$ if and only if
$G[V_{\mathbb{R}}(G)]=G(\mathbb{R})=G(E')$ and $g(G')=0$.
\renewcommand{\qedsymbol}{ }
\end{proof}

Now we are going to prove $s(G) \equiv g(G)+1 \pmod{2}$. We write
$\overline{V}_{\mathbb{R}}(G)$ to denote the set of vertices of
$G(E')$. We know $\overline{V}_{\mathbb{R}}(G)$ can be considered as
a quotient set of $V_{\mathbb{R}}(G)$. Clearly one has

\begin{equation}\label{formule6}
    g(G(E'))=g(G[V_{\mathbb{R}}(G)])+(\sharp
    (V_{\mathbb{R}}(G))-\sharp (\overline{V}_{\mathbb{R}}(G))).
\end{equation}

\noindent Because of Formula \ref{formule3} this implies

\begin{equation}\label{formule7}
    g(G(E')) \equiv s(G)+1+(\sharp (V_{\mathbb{R}}(G))-\sharp
    (\overline{V}_{\mathbb{R}}(G))) \pmod{2}.
\end{equation}

We make a partition $V_{\mathbb{R}}(G)=V_{\mathbb{R}}(G)_0\cup
V_{\mathbb{R}}(G)_1$ such that a real vertex $v$ belongs to
$V_{\mathbb{R}}(G)_0$ if and only if for each edge $e$ with $v\in
\psi (e)$ one has $\psi (e)\subset V_{\mathbb{R}}(G)$. This is
equivalent to $\{ v\}$ being a component of $G'=G\setminus E'$. We
have an induced partition
$\overline{V}_{\mathbb{R}}(G)=V_{\mathbb{R}}(G)_0\cup
\overline{V}_{\mathbb{R}}(G)_1$ (here
$\overline{V}_{\mathbb{R}}(G)_1$ is the quotient set of
$V_{\mathbb{R}}(G)_1$ using the identification mentioned before) and
we obtain

\begin{equation}\label{formule8}
    g(G(E')) \equiv s(G)+1+(\sharp (V_{\mathbb{R}}(G)_1)-\sharp
    (\overline{V}_{\mathbb{R}}(G)_1)) \pmod{2}.
\end{equation}

\begin{figure}[h]
\begin{center}
\includegraphics[height=3 cm]{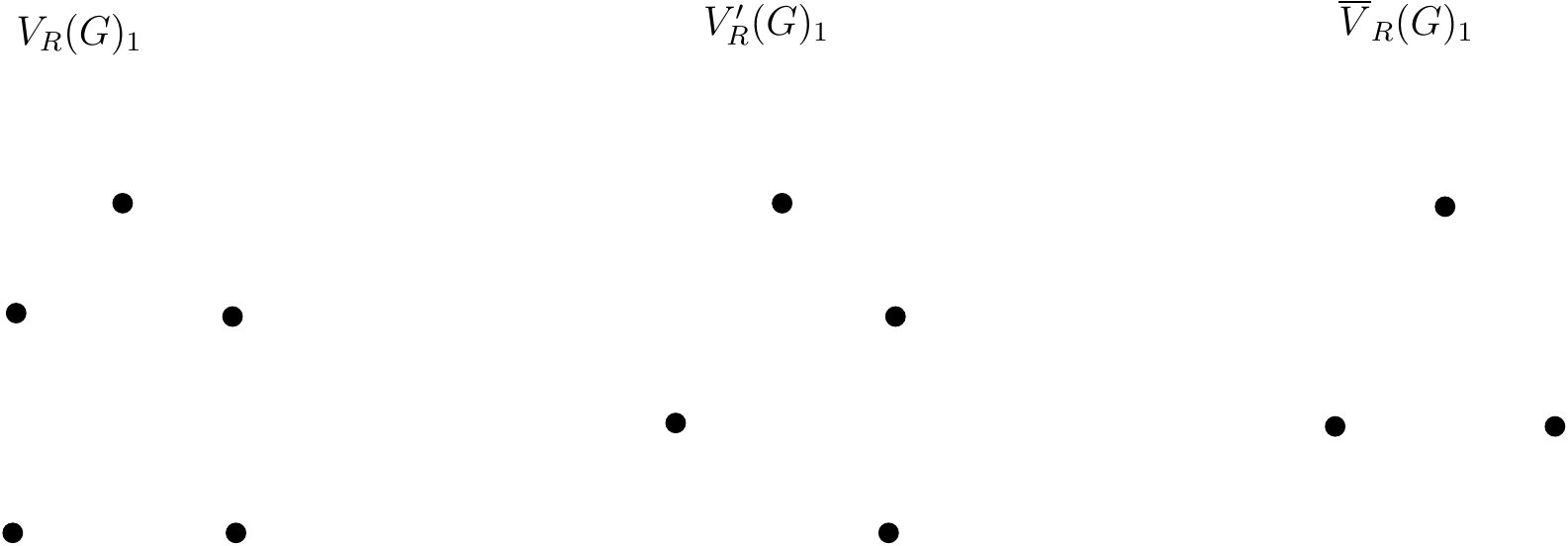}
\caption{Identification of real vertices of the example.}
\end{center}
\end{figure}

Let $S'$ be the set of edges $e$ of $G'$ satisfying $\psi (e)$
contains a real vertex. This set is invariant under $\iota _E$ and
it contains no real edge, hence $\sharp (S')$ is even. We consider
the decomposition

\begin{equation}\label{formule9}
    g(G')=g(G'(S'))+g(G'\setminus S').
\end{equation}

Let $V'_{\mathbb{R}}(G)_1$ be the quotient set of
$V_{\mathbb{R}}(G)_1$ obtained by identifying vertices $v_1$ and
$v_2$ from $V_{\mathbb{R}}(G)_1$ if and only if they are ends of a
walk with edges belonging to $S'$. There is a natural bijection
between the set of connected components of the graph $G'[S']$ and
$V'_{\mathbb{R}}(G)_1$. Moreover $\overline{V}_{\mathbb{R}}(G)_1$
can be considered as a quotient set of $V'_{\mathbb{R}}(G)_1$. Also
$\sharp (V(G'[S'])\setminus V_{\mathbb{R}}(G)_1)$ is even. Those
remarks imply

\begin{equation}\label{formule10}
    g(G'[S']) \equiv \sharp (V'_{\mathbb{R}}(G)_1)- \sharp (V
    _{\mathbb{R}}(G)_1) \pmod{2}.
\end{equation}

The graph $G'(S')$ is obtained from $G'[S']$ by identifying vertices
$v_1$ and $v_2$ of $G'[S']$ if and only if there is a walk $\Gamma$
in $G'\setminus S'$ with ends $v_1$ and $v_2$. Since $\psi (e)$
contains no real vertex for an edge $e$ of $G'\setminus S'$ it
follows $v_1$ and $v_2$ are non-real vertices. Since $a(G)=0$ also
$v_2\neq \overline{v_1}$. On the other hand, $\overline{\Gamma}$ is
a walk in $G'\setminus S'$ with ends $\overline{v_1}$ and
$\overline{v_2}$. Hence $\overline{v_1}$ and $\overline{v_2}$ as
vertices of $G'[S']$ are also identified in $G'(S')$. This proves
$\sharp (V(G'[S']))-\sharp (V(G'(S')))$ is even.

\begin{figure}[h]
\begin{center}
\includegraphics[height=4 cm]{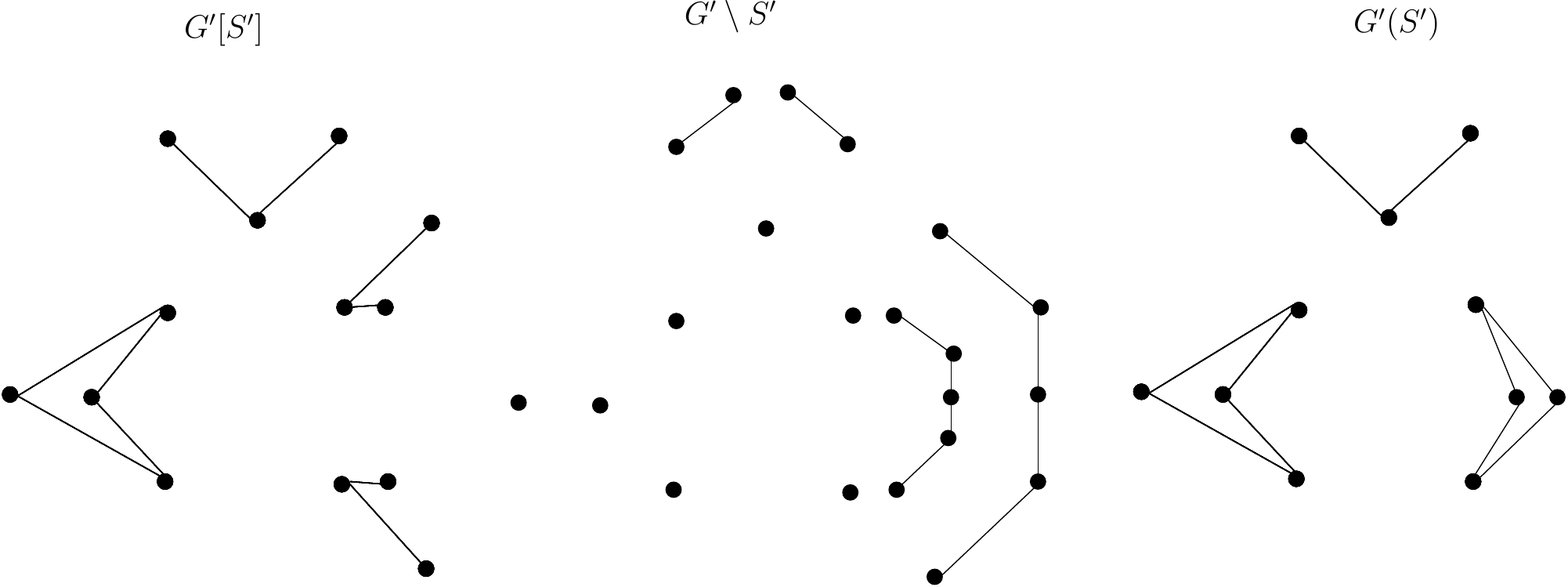}
\caption{Graph $G'(S')$ of the example.}
\end{center}
\end{figure}

Assume $C_1$ and $C_2$ are two different connected components of
$G'[S']$ represented by $v'_1$ and $v'_2$ in $V'_{\mathbb{R}}(G)_1$.
Assume $C_i$ has vertex $v_i$ such that $v_1$ and $v_2$ are
identified in $G'(S')$. This is equivalent to the existence of a
walk $\Gamma$ in $G'\setminus S'$ that connects $v_1$ to $v_2$.
Hence it is equivalent to the fact that $v'_1$ and $v'_2$ do define
the same points in $\overline{V}_{\mathbb{R}}(G)_1$. This implies
that the difference between the number of connected components of
$G'(S')$ and $G'[S']$ is equal to $\sharp
(V'_{\mathbb{R}}(G)_1)-\sharp (\overline{V}_{\mathbb{R}}(G)_1)$.
Those arguments imply

\begin{equation}\label{formule11}
    g(G'(S'))-g(G'[S']) \equiv \sharp (V'_{\mathbb{R}}(G)_1)-\sharp
    (\overline{V}_{\mathbb{R}}(G)_1) \pmod{2}.
\end{equation}

\noindent Combining this with Formula \ref{formule10} we obtain

\begin{equation}\label{formule12}
    g(G(S')) \equiv \sharp (V_{\mathbb{R}}(G)_1)-\sharp
    (\overline{V}_{\mathbb{R}}(G)_1) \pmod{2}.
\end{equation}

By construction $G'\setminus S'$ is a real subgraph of $G$. Let $H$
be a connected component of $G'\setminus S'$ different from a vertex
then $\overline{H}$ is also a connected component of $G'\setminus
S'$. Since $G'\setminus S'$ contains no real edge and each real
vertex is an isolated component of $G'\setminus S'$ we know $H$
contains no real vertex and no real edge. Since $a(G)=0$ it follows
$H\neq \overline{H}$. But $H$ is isomorphic to $\overline{H}$ hence
$g(H)=g(\overline{H})$. This implies $g(G'\setminus S')$ is even,
therefore formulas \ref{formule9} and \ref{formule12} imply

\begin{equation}\label{formule13}
    g(G') \equiv \sharp (V_{\mathbb{R}}(G)_1)-\sharp
    (\overline{V}_{\mathbb{R}}(G)_1) \pmod{2}.
\end{equation}

\noindent Combining this with Formulas \ref{formule4} and
\ref{formule8} we obtain the desired formula $g(G) \equiv s(G)+1
\pmod{2}$.

\begin{proof}[Step 5]
Induction argument on the number of connected components of
$G(\mathbb{R})$.
\renewcommand{\qedsymbol}{ }
\end{proof}

Now we assume $G(\mathbb{R})$ is not connected and fix a connected
component $G_1$ of $G(\mathbb{R})$. Since $G$ is connected there
exists another connected component $G_2$ of $G(\mathbb{R})$ such
that there is a walk $\Gamma$ in $G$ that connects a vertex $v_1$ of
$G_1$ to a vertex $v_2$ of $G_2$ having no real inner vertex. We
also can assume that all vertices of $\Gamma$ are mutually
different. Then $\overline{\Gamma}$ satisfies the same properties.
Moreover $\overline{\Gamma}$ is different from $\Gamma$, indeed, the
only edge $e$ of $\Gamma$ containing $v_1$ is non-real and
$\overline{e}$ is the only edge of $\overline{\Gamma}$ containing
$v_1$. Let $S=\Gamma \cup \overline{\Gamma}$. Clearly $S$ is a real
subgraph of $G$, hence also $G\setminus E(S)$ is a real subgraph of
$G$. Moreover $(G\setminus E(S))(\mathbb{R})=G(\mathbb{R})$. Let
$G'_1; \cdots ; G'_t$ be the connected subgraphs of $G\setminus
E(S)$. If $G'_i$ contains some real vertex of $G$ then $G'_i$ is a
real subgraph of $G$ and $a(G)=0$ implies $a(G'_i)=0$. Assume $1\leq
x\leq t$ such that $G'_i$ contains a real vertex of $G$ if and only
if $1\leq i\leq x$. Also assume $v_1$ is a vertex of $G'_1$ and in
case $v_2$ is not a vertex of $G'_1$, $v_2$ is a vertex of $G'_2$.
We are going to assume that the theorem holds for those real
subgraphs $G'_1; \cdots; G'_x$ and we are going to prove that this
implies the theorem holds for $G$. In case $x\geq 2$ for each $1\leq
i\leq x$ the graph $G'_i(\mathbb{R})$ has less connected components
than $G(\mathbb{R})$ and we can use the induction hypothesis. In
case $x=1$, $G'_1$ is a proper real subgraph of $G$ such that
$G'_1(\mathbb{R})=G(\mathbb{R})$. In particular $G'_1(\mathbb{R})$
is not connected and we can apply the same construction to $G'_1$
and continuing in this way we arrive at a situation such that $x\geq
2$. Again we can use the induction hypothesis.

\begin{figure}[h]
\begin{center}
\includegraphics[height=6 cm]{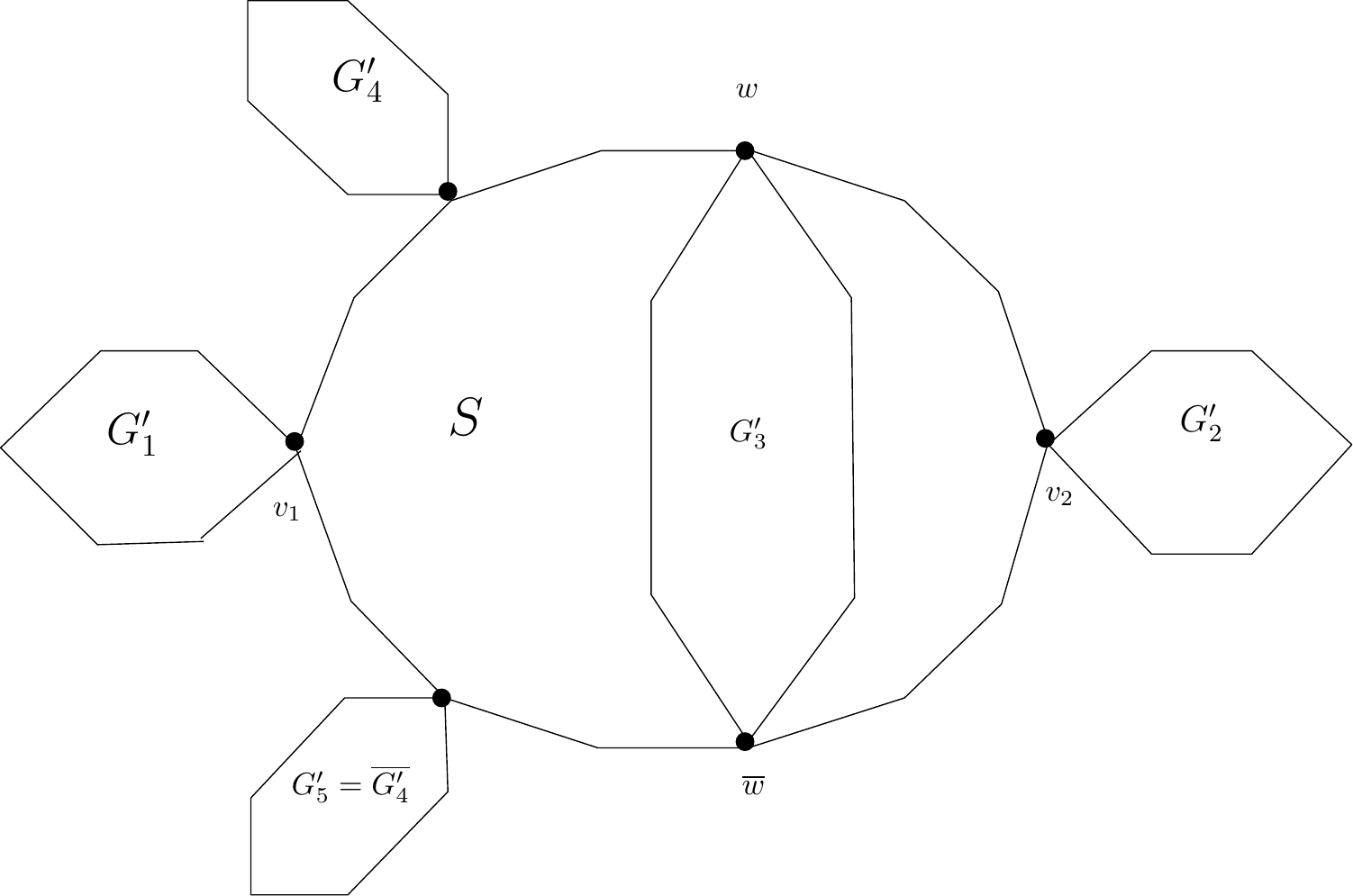}
\caption{Connected components of $G\setminus E(S).$}
\end{center}
\end{figure}

For $1\leq i\leq t$ let $W_i=V(S)\cap V(G'_i)$. Since $G$ is
connected $W_i\neq \emptyset$ for $1\leq i\leq t$ and since $G'_1;
\cdots ; G'_t$ are different connected components of $G\setminus S$
one has $W_i\cap W_j=\emptyset$ for $i\neq j$. Moreover $G(E(S))$ is
obtained from $S=G[E(S)]$ by identifying vertices $w_1$ and $w_2$ of
$S$ if and only if there exists $1\leq i\leq t$ such that $w_1;
w_2\in W_i$. This implies

\begin{equation}\label{formule14}
    g(G(E(S)))=g(S)+\left( \sum_{i=1}^t \sharp (W_i)\right) -t.
\end{equation}

\noindent We use the decomposition

\begin{equation}\label{formula15}
    g=g(G)=g(G(E(S)))+g(G\setminus E(S))=g(G(E(S)))+\sum_{i=1}^t
    g(G'_i).
\end{equation}

\noindent For $1\leq i\leq x$ we write $s_i=s(G'_i)$ and we assume
$s_i \equiv g(G'_i)+1 \pmod{2}$ and $s_i\leq g(G'_i)+1$ for $1\leq
i\leq x$. This implies

\begin{equation}\label{formule16}
    g\geq g(G(E(S)))+\sum_{i=1}^x (s_i-1)+\sum_{i=x+1}^t g(G'_i).
\end{equation}

\noindent with equality if and only if $s_i=g(G_i)+1$ for $1\leq
i\leq x$. Also by definition $s(G)=\sum_{i=1}^x s_i$.

Let $1\leq i\leq x$ and assume $w\in W_i\setminus \{ v_1;v_2 \}$.
Since $\overline{G_i}=G_i$ it follows $\overline{w}\in W_i$. In case
$v_2 \in V(G'_1)$ it implies $\sharp (W_i)$ is even for $1\leq i\leq
x$ and in case $v_2\in V(G'_2)$, $\sharp (W_1)$ and $\sharp (W_2)$
are odd and $\sharp (W_i)$ is even for $3\leq i\leq x$. In both
cases we obtain $\sum_{i=1}^x \sharp (W_i)$ is even. In case $v_2
\in V(G'_1)$ it also implies $\sum_{i=1}^x \sharp (W_i)\geq 2x$ and
in case $v_2\in V(G'_2)$ it implies $\sum_{i=1}^x \sharp (W_i)\geq
2x-2$. Since $g(G[S])\geq 1$ it follows, combining formulas
\ref{formule14} and \ref{formule16} that

\begin{equation}\label{formule17}
    g(G)\geq 1+2x-x+\sum_{i=1}^x (s_i-1)=s(G)+1
\end{equation}

\noindent in case $v_2\in V(G'_1)$ and

\begin{equation}\label{formule18}
    g(G)\geq 1+2x-2-x+\sum_{i=1}^x (s_i-1)=s(G)-1
\end{equation}

\noindent in case $v_2\in V(G'_2)$. We obtain $s(G)\leq g(G)+1$ and
moreover the computation implies the following claim we are going to
use in the proof of Theorem \ref{theorem3}.

\begin{proof}[Claim]
If $s(G)=g(G)+1$ then $v_2\in V(G_2)$, $\sharp (W_1)=\sharp
(W_2)=1$, $\sharp (W_i)=2$ for $3\leq i\leq x$, $\sharp (W_i)=1$ for
$x+1\leq i\leq t$, $g(G'_i)=s(G'_i)-1$ for $1\leq i\leq x$ (in
particular $a(G'_i)=0$), $g(G'_i)=0$ for $x+1\leq i\leq t$ and
$g(S)=1$ (hence $S$ is a cycle).
\renewcommand{\qedsymbol}{ }
\end{proof}

In order to finish the proof of the theorem we need to prove $s(G)
\equiv g(G)+1 \pmod{2}$. If $S$ would not be a cycle then there
would be an inner vertex $w$ of $\Gamma$ such that $\overline{w}$ is
also a vertex of $\Gamma$. The subwalk of $\Gamma$ that connects $w$
to $\overline{w}$ would imply $a(G)=0$, hence a contradiction. By
induction we have $s_i \equiv g(G'_i)+1 \pmod{2}$ for $1\leq i\leq
x$. Since $G'_i$ contains no real vertex for $x+1\leq i\leq t$ also
$\overline{G'_i}\neq G'_i$ in that case. This implies
$\sum_{i=x+1}^t g(G'_i) \equiv 0 \pmod{2}$, $t \equiv x \pmod{2}$
and $\sum_{i=x+1}^t \sharp (W_i) \equiv 0 \pmod{2}$. We already
noticed that $\sum_{i=1}^x \sharp (W_i) \equiv 0 \pmod{2}$.
Therefore, from formulas \ref{formule14} and \ref{formula15} it
follows that $g(G) \equiv 1-x+\sum_{i=1}^x (s_i-1) \pmod{2}$, hence
$g(G) \equiv s(G)+1 \pmod{2}$.

\end{proof}

\begin{example}\label{voorbeeld1}
\normalfont Let $g$, $a$ and $s$ be nonnegative integers such that
$0\leq s\leq g+1$ with $s \equiv g+1 \pmod{2}$ and $a\in \{ 0;1\}$
such that $a=0$ if $s=g+1$ and $a=1$ if $s=0$ (hence $g$ is odd in
this case). We give an example of a graph $G$ with a real structure
satisfying $g(G)=g$, $a(g)=a$ and $s(G)=s$.

First assume $a=0$ (hence $s\geq 1$) and write $g+1-s=2x$ (with
$x\geq 0$ because $s\leq g+1$). Take vertices $v_1; \cdots ; v_s$
and two different edges $e_i; \overline{e_i}$ for $1\leq i\leq s-1$
such that $\psi (e_i)=\psi (\overline{e_i})=\{ v_i; v_{i+1} \}$.
Then take vertices $w_1; \cdots ; w_x; \overline{w_1}; \cdots ;
\overline{w_x}$ and edges $f_i; f'_i; \overline{f_i};
\overline{f'_i}$ ($1\leq i\leq x$) with $\psi (f_1)=\psi (f'_1)=\{
v_1; w_1\}$, $\psi (\overline{f_1})=\psi (\overline{f'_1})=\{ v_1;
\overline{w_1}\}$, $\psi (f_i)=\psi (f'_i)=\{ w_{i-1};w_i\}$ and
$\psi (\overline{f_i})=\psi (\overline{f'_i})=\{ \overline{w_{i-1}};
\overline{w_i} \}$. We obtain a graph $G$ such that $g(G)=s-1+2x=g$.
Define the real structure on $G$ such that $^{\overline{ \text{ }}}$
is conjugation and with $\overline{v_i}=v_i$, then $G(\mathbb{R})$
is the disjoint union of the trivial graphs $\{ v_1 \}\cup \{ v_2 \}
\cup \cdots \cup \{ v_s\}$, hence $s(G)=s$. Clearly $a(G)=0$.

\begin{figure}[h]
\begin{center}
\includegraphics[height=4 cm]{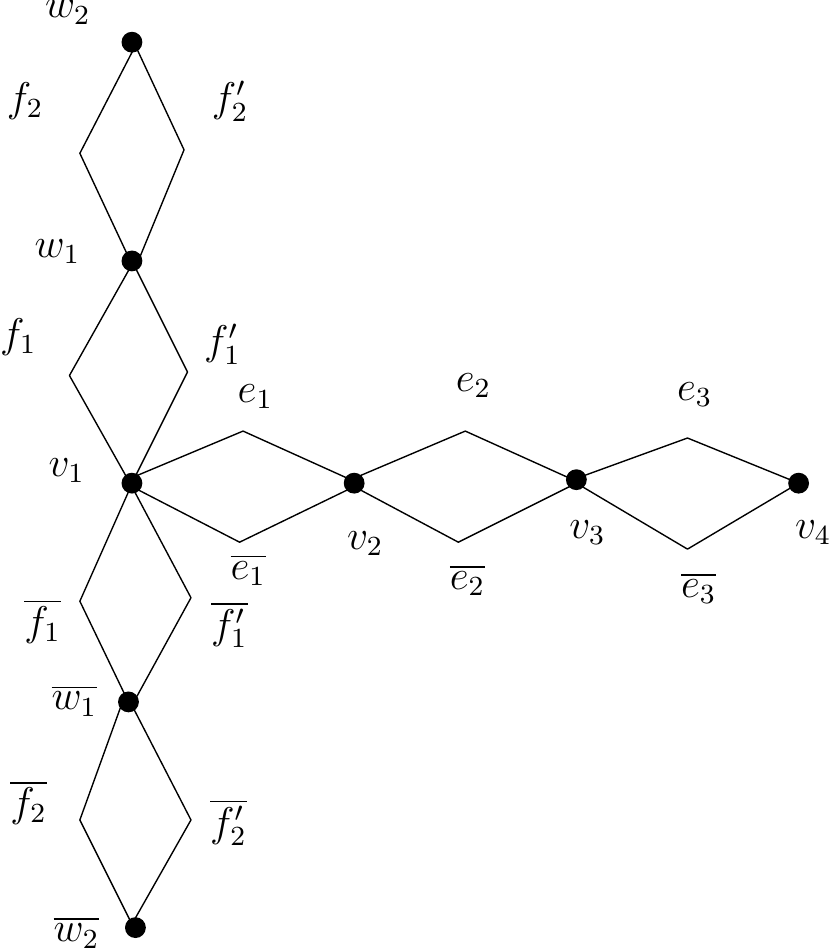}
\caption{Graph with $a=0$.}
\end{center}
\end{figure}

Next assume $s=0$, hence $a=1$ and $g$ is odd. Write $g+1=2x$ and
take two vertices $v; \overline{v}$ and $2x$ edges $e_1; \cdots ;
e_x; \overline{e_1} ; \cdots ; \overline{e_x}$ with $\psi (e_i)=\psi
(\overline{e_i})=\{ v;\overline{v}\}$ for $1\leq i\leq x$. Clearly
$g(G)=2x-1=g$. Defining the real structure on $G$ such that
$^{\overline{ \text{ }}}$ denotes conjugation we obtain $s(G)=0$ and
clearly $a(G)=1$.

Finally take $s\neq 0$ but $a=1$. Write $g+1-s=2x$ with $x\geq 1$
because $s\leq g-1$. Take $s$ vertices $v_1; \cdots ; v_s$ and
$2s-2$ edges $e_1; \overline{e_1} ; \cdots ; e_{s-1} ;
\overline{e_{s-1}}$ with $\psi (e_i)=\psi (\overline{e_i})=\{ v_i;
v_{i+1} \}$ for $1\leq i\leq s-1$. Take two more vertices $v$ and
$\overline{v}$ and $2x+2$ more edges $f; \overline{f}; f_1;
\overline{f_1}; \cdots ; f_x; \overline{f_x}$ with $\psi (f)=\{ v_s
; v\}$, $\psi (\overline{f})=\{ v_s; \overline{v}\}$ and $\psi
(f_i)=\psi (\overline{f_i})=\{ v; \overline{v}\}$ for $1\leq i\leq
x$. Then $g(G)=s-1+1+2x-1=g$. Take a real structure on $G$ such that
$^{\overline{ \text{ }}}$ denotes conjugation and
$\overline{v_i}=v_i$ for $1\leq i\leq s$. As in the first part of
the example we have $s(G)=s$ and since $x\geq 1$, by construction we
also have $a(G)=0$.

\begin{figure}[h]
\begin{center}
\includegraphics[height=4 cm]{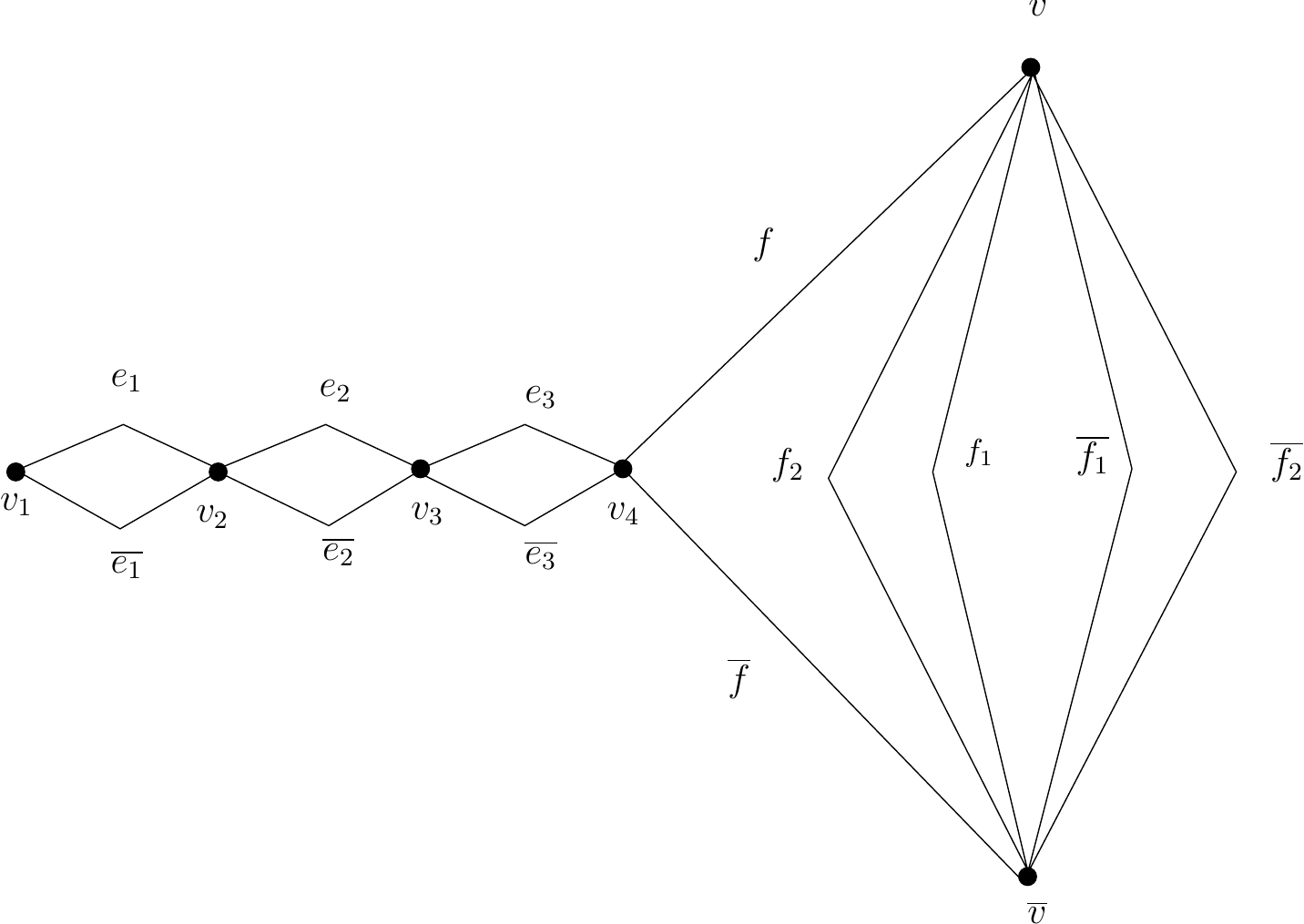}
\caption{Graph with $a=1$.}
\end{center}
\end{figure}

\end{example}

\section{Real linear systems}\label{section3}

We start by recalling some definitions concerning divisors and
linear systems on graphs as described in e.g. \cite{ref1}. A
\emph{divisor} on a graph $G$ is a formal integer linear combination
of vertices of $G$. In case $D$ is such a divisor and $v$ is a
vertex we write $D(v)$ to denote the coefficient of $D$ at $v$,
hence $D=\sum_{v\in V(G)}D(v)v$. The \emph{support} $\Supp (D)$ of
$D$ is the subset of $V(G)$ of vertices $v$ satisfying $D(v)\neq 0$.
A divisor is called \emph{effective} if $D(v)\geq 0$ for each vertex
$v\in V(G)$ and in that case we write $D\geq 0$. For a divisor $D$
we call $\sum_{v\in V(G)}D(v)$ the degree $\deg (D)$ of $D$. We
write $\Div (G)$ to denote the group of divisors on $G$ (hence it is
the free abelian group on $V(G)$). We write $\Div ^0(G)$ to denote
the subgroup of divisors on $G$ of degree 0.

Let $\mathcal{M}(G)$ be the abelian group of integer-valued
functions $f$ on $V(G)$. For $f\in \mathcal{M}(G)$ we define the
\emph{principal divisor} $\Delta (f)$ as follows. For $v\in V(G)$ we
define

\begin{equation}\label{formule19}
    \Delta (f)(v)=\sum_{ \substack{ e\in E(G)\\ \psi (e)=\{ v;
    w\}}}(f(w)-f(v)).
\end{equation}

\noindent Clearly $\Delta (f)\in \Div ^0(G)$ and $\{ \Delta (f) :
f\in \mathcal{M}(G)\}$ is a subgroup $\Prin (G)$ of $\Div ^0(G)$. We
say two divisors $D_1; D_2\in \Div (G)$ are linearly equivalent if
and only if $D_2 -D_1\in \Prin (G)$. In particular $\deg (D_1)=\deg
(D_2)$ and we write $D_1 \sim D_2$. For a divisor $D\in \Div(G)$ we
define the \emph{complete linear system} $| D|$ as follows:

\begin{displaymath}
    | D| = \{ D' \in \Div (G) : D'\geq 0 \text { and } D'\sim D \}
\end{displaymath}

\noindent We also introduce a notion for the \emph{rank} $\rk (D)$
of $D$ as follows. In case $| D|=\emptyset$ we put $\rk (D)=-1$. So
assume $| D|\neq \emptyset$, then $\rk (D)$ is the maximal integer
value $r$ such that for each $E\geq 0$ with $\deg (E)=r$ there
exists $D'\in | D|$ with $D'-E\geq 0$. For a divisor $D$ on a smooth
curve $X$ the corresponding notion is the same as the dimension of
the complete linear system $| D|$.

Now assume $G$ has a real structure. Associated to $D=\sum_{v\in
V(G)}D(v)v$ there is a conjugated divisor $\overline{D}=\sum_{v\in
V(G)}D(\overline{v})v$, hence $\overline{D}(v)=D(\overline{v})$. We
say $D$ is a \emph{real divisor} if $D=\overline{D}$. In case $f\in
\mathcal{M}(G)$ we define $\overline{f}\in \mathcal{M}(G)$ by
$\overline{f}(v)=f(\overline{v})$. Clearly
$\Delta(\overline{f})=\overline{\Delta(f)}$.

\begin{definition}\label{definition7}
Let $D$ be a real divisor on a graph $G$ with a real structure. The
\emph{real rank} $\rk_{\mathbb{R}} (D)$ is the maximal integer value
$r$ such that for each real divisor $E\geq 0$ with $\deg (E)=r$
there exists a real divisor $D'\in | D|$ with $D'-E\geq 0$.
\end{definition}

In the case of a real curve $X$ the real dimension of a complete
linear system associated to a real divisor $D$ is equal to the
dimension of the complete linear system $| D|$ on the complex curve
$X_{\mathbb{C}}$ (see e.g. \cite{ref5}*{p. 200}). In the graph-case
this becomes only one inequality.

\begin{proposition}\label{propositie1}
Let $G$ be a graph with a real structure and let $D$ be a real
divisor on $G$. One has $\rk _{\mathbb{R}}(D)\geq \rk (D)$.
\end{proposition}

\begin{proof}
In case $\rk (D)=-1$ there is nothing to prove, so assume $\rk
(D)\geq 0$.

Let $E$ be a real effective divisor of degree $r$ on $G$. There
exists $D''\in | D |$ such that $D''-E\geq 0$. Let $f\in \mathcal{M}
(G)$ with $D+\Delta (f)=D''$. Since $\overline{D}=D$ we obtain
$D+\Delta (\overline{f})=\overline{D''}$, hence $\overline{D''}\in
|D|$.

Let $g=\max \{ f; \overline{f}\}$ and let $D'=D+\Delta (g)$. We are
going to prove that $D'\geq 0$ (hence $D'\in | D|$);
$\overline{D'}=D'$ (hence $D'$ is real) and $D'\geq E$. This will
imply the proposition.

Since $D''\geq 0$, for each $v\in V(G)$ one has

\begin{equation}\label{formule20}
    D(v)+\sum_{ \substack{ e\in E(G) \\ \psi (e)=\{ v;w
    \}}}(f(w)-f(v))=D''(v)\geq 0.
\end{equation}

\noindent Using $\overline{v}$ instead of $v$ and taking into
account that $D(\overline{v})=D(v)$, we obtain

\begin{equation}\label{formule21}
    D(v)+\sum_{ \substack{ e\in E(G) \\ \psi (e)=\{ v;w \}}}
    (\overline{f}(w)-\overline{f}(v))=D''(\overline{v})\geq 0
\end{equation}

Assume $f(v)\geq f(\overline{v})$, hence $g(v)=f(v)$. For each $e\in
E(G)$ with $\psi (e)=\{ v;w\}$ one has $g(w)-g(v)\geq f(w)-f(v)$
therefore Formula \ref{formule20} implies $D(v)+\Delta
(g)(v)=D'(v)\geq D''(v)\geq 0$. Assume $\overline{f}(v)\geq f(v)$
hence $g(v)=\overline{f}(v)$. Now for $e$ as before we have
$g(w)-g(v)\geq \overline{f}(w)-\overline{f}(v)$. Again, using
Formula \ref{formule21} we obtain $D(v)+\Delta (g)(v)=D'(v)\geq
D''(\overline{v})\geq 0$. This proves $D'\in | D|$.

By definition $\overline{g}=g$, hence $\overline{\Delta (g)}=\Delta
(g)$ and therefore $\overline{D'}=\overline{D+\Delta (g)}=D'$. This
proves D' is a real divisor.

For each vertex $v$ we have $D''(v)\geq E(v)$. In case $v$ is a real
vertex $g(v)=f(v)$ and in such case we obtained $D'(v)\geq D''(v)$
hence $D'(v)\geq E(v)$. Assume $v$ is not a real vertex. Note that
$D''(\overline{v})\geq E(\overline{v})=E(v)$. We can assume
$g(v)=f(v)$ (otherwise we use $\overline{v}$ instead of $v$) hence
$D'(v)\geq D''(v)$ and therefore $D'(v)\geq E(v)$. Since $D'$ and
$E$ are both real divisors we also have $D'(\overline{v})\geq
E(\overline{v})$.

\end{proof}

\begin{example}\label{voorbeeld2}
\normalfont Let $G'$ be a graph with $g(G')\geq 1$. Fix $v'\in
V(G')$. We construct a new graph $G$ by taking two copies $G'_1$ and
$G'_2$ of $G'$ (hence $v'_i\in V(G'_i)$ is the corresponding copy of
$v'$), a new vertex $v$ and two new edges $e_1; e_2$ with $\psi
(e_i)=\{ v; v'_i\}$. On $G$ we define a real structure such that
$\iota$ induces the identification of $G'_1$ and $G'_2$ using their
identification with $G'$, $\overline{e_1}=e_2$ and $v$ is a real
vertex. Then $D=v$ is a real divisor and clearly $\rk_{\mathbb{R}} (
D)=1$. Since $g(G)\neq 0$ there is no linear system of degree 1 and
dimension at least 1, hence $\rk (D)=0$ (see e.g. \cite{ref6}*{Lemma
1.1}). This proves that in general Proposition \ref{propositie1}
cannot have equality.

\end{example}

\begin{lemma}\label{lemma1}
Assume $D_1$ and $D_2$ are real divisors on a graph $G$ with a real
structure and $f\in \mathcal{M} (G)$ with $\Delta (f)=D_1-D_2$, then
$\overline{f}=f$
\end{lemma}

\begin{proof}
Since $\Delta (f)=D_2-D_1$ and $D_1$ and $D_2$ are real one has
$\overline{\Delta (f)}=\Delta (\overline{f})=D_2-D_1$. This implies
$\Delta (f-\overline{f})=0$. However, for $g\in \mathcal{M} (G)$ the
equation $\Delta (g)=0$ implies $g$ is a constant function on
$V(G)$. If not, let $V_0\subset V(G)$ be the subset of $V(G)$
consisting of all $v\in V(G)$ such that $g$ has minimal value at
$v$. Since $V_0\neq V(G)$ and $G$ is connected there exists $v_0\in
V_0$, $v\in V(G)\setminus V_0$ and $e\in E(G)$ with $\psi (e)=\{
v_0;v\}$. Of course $\Delta (g)(v_0)>0$ for such vertex $v_0$,
contradicting $\Delta (g)=0$.

It follows $f-\overline{f}$ is a constant function on $V(G)$. For
each vertex $v$ of $V(G)$ one has
$(f-\overline{f})(v)=(f-\overline{f})(\overline{v})$ hence

\begin{displaymath}
f(v)-\overline{f}(v)=f(\overline{v})-\overline{f}(\overline{v})=f(\overline{v})-f(v)
\end{displaymath}

\noindent hence $2f(v)=2f(\overline{v})$. This implies
$f(v)=f(\overline{v})=\overline{f}(v)$ and therefore
$f=\overline{f}$.
\end{proof}

In case $X$ is a real curve and $C$ is a connected component of
$X(\mathbb{R})$ the parity of the degree of the restriction of real
divisors on $X$ to $C$ is invariant under linear equivalence. This
is a basic fact in the study of real linear systems on real curves
and in the study of their real projective embeddings. The next
theorem shows this also holds for graphs with a real structure.

\begin{theorem}\label{theorem2}
Let $G$ be a graph with a real structure and let $D_1$ and $D_2$ be
two linearly equivalent real divisors on $G$. Let $G'$ be a
connected component of $G(\mathbb{R})$, then $\deg (D_1| _{G'})
\equiv \deg (D_2| _{G'}) \pmod{2}$.
\end{theorem}

\begin{proof}
Because of Lemma \ref{lemma1} there exists $f\in \mathcal{M} (G)$
with $\overline{f}=f$ such that $D_2=D_1+\Delta (f)$. Let $f'$ be
the restriction of $f$ to $V(G')$. Let $\Delta (f')$ be the
associated principal divisor on $G'$, then $\deg (D_1|_{G'})=\deg
(D_1|_{G'}+\Delta (f'))$. Let $D'_2=D_1|_{G'}+\Delta (f')$. For
$v\in V(G')$ one has

\begin{displaymath}
    D_2(v)=D_1(v)+\sum_{ \substack{ e\in E(G') \\ \psi (e)=\{
    v;w\}}} (f'(w)-f'(v))+\sum_{ \substack{ e\in E(G)\setminus E(G')
    \\ \psi (e)=\{ v;w\}}} (f(w)-f(v))
\end{displaymath}

\noindent hence

\begin{displaymath}
    D_2(v)=D'_2(v)+\sum_{ \substack{ e\in E(G)\setminus E(G')
    \\ \psi (e)=\{ v;w\}}} (f(w)-f(v)).
\end{displaymath}

\noindent For each $e\in E(G)\setminus E(G')$ with $\psi (e)=\{
v;w\}$ one has $\overline{e}\neq e$ and $\psi (\overline{e})=\{v;
\overline{w}\}$ and
$f(\overline{w})-f(v)=\overline{f}(w)-f(v)=f(w)-f(v)$. This proves
$\sum_{ \substack{ e\in E(G)\setminus E(G') \\ \psi (e)=\{ v;w\}}}
(f(w)-f(v))$ is even, hence $D_2(v) \equiv D'_2(v) \pmod{2}$. So we
obtain

\begin{displaymath}
    \deg (D_2|_{G'})=\sum_{v\in V(G')}D_2(v) \equiv \sum_{v\in
    V(G')} D'_2(v) \pmod{2}
\end{displaymath}

\noindent Since $\sum_{v\in V(G')}D'_2(v)=\deg (D'_2)=\deg
(D_1|_{G'})$, this finishes the proof.

\end{proof}

In the case of smooth real curves $X$ the parity of the canonical
linear system on $X$ (which is a real linear system) is even on each
connected component of $X(\mathbb{R})$ (see e.g.
\cite{ref3}*{Corollary 4.3}). In the case of graphs there is a
distinguished canonical divisor

\begin{displaymath}
    K_G=\sum_{v\in V(G)} (\val v -2)v
\end{displaymath}

\noindent (here $\val v$ is the number of edges $e\in E(\mathbb{R})$
with $v\in \psi (e)$).

\begin{proposition}\label{propositie2}
Let $G$ be a graph with a real structure and let $G'$ be a connected
component of $G(\mathbb{R})$, then $\deg (K_G|_{G'})$ is even.
\end{proposition}

\begin{proof}
On the graph $G'$ there is a canonical divisor

\begin{displaymath}
    K_{G'}=\sum_{v\in V(G')} (\deg _{G'}(v)-2)v
\end{displaymath}

\noindent One has $\deg (K_G')=2g(G')-2$, hence it is even. For
$v\in V(G')$ the difference $K_G(v)-K_{G'}(v)$ is equal to the
number of non-real edges $e\in E(G)$ with $v\in \psi (e)$. Of course
this number of edges is even, hence $K_G(v)-K_{G'}(v)$ is even. This
proves $\deg (K_G|_{G'})$ is even.

\end{proof}

\begin{definition}\label{definitie8}
A graph $G$ with a real structure is called an \emph{M-graph} if $G$
has no isolated real edge and $s(G)=g(G)+1$.
\end{definition}

\begin{remark}
If $G$ satisfies $s(G)=g(G)+1$ but $G$ has isolated real edges then
there is an associated M-graph obtained as explained in step 1 of
the proof of Theorem \ref{theorem1}
\end{remark}

\begin{definition}\label{definition9}
A real divisor $D$ on a graph $G$ with a real structure is called
\emph{totally real} if $\Supp (D)\subset V_{\mathbb{R}}(G)$.
\end{definition}

The following result concerning real divisors on M-graphs is a very
strong one and the corresponding statement for real divisors on
M-curves does not hold.

\begin{theorem}\label{theorem3}
Let $G$ be an M-graph and let $D$ be an effective real divisor on
$G$ then $D$ is linearly equivalent to a totally real effective
divisor on $G$.
\end{theorem}

\begin{proof}

We are going to prove a stronger claim: if $v+\overline{v}$ is a
non-real vertex pair on $G$ then there is a real vertex $w$ such
that $v+\overline{v}\sim 2w$. The proof is going to make use of
notations and claims from the proof of Theorem \ref{theorem1}.

Of course, in case $G=G(\mathbb{R})$ there is nothing to prove, so
we assume $G\neq G(\mathbb{R})$. First assume $G(\mathbb{R})$ is
connected. As mentioned in the claim inside step 4 of the proof of
Theorem \ref{theorem1} we obtained
$G(E')=G[V_{\mathbb{R}}(G)]=G(\mathbb{R})$, hence $g(G(E'))=s(G)-1$,
and $g(G')=0$. In particular each edge $e\in E(G)$ satisfying $\psi
(e)\subset V_{\mathbb{R}}(G)$ is a real edge. Let $v$ be a non-real
vertex of $G$. There is a walk $\Gamma$ that connects $v$ to a real
vertex $v'$ such that no inner vertex of $\Gamma$ is a real vertex.
Let $V(v)$ be the set of vertices $w$ of $G$ such that there is a
walk with ends $v$ and $w$ and containing no real vertex. Let $T(v)$
be subgraph of the induced subgraph $G[V(v)\cup \{ v'\}]$ obtained
by omitting the loops at $v'$. It is a subgraph of $G'$. Since
$g(G')=0$ the graph $T(v)$ is a tree. Since
$G(E')=G[V_{\mathbb{R}}(G)]$ there is no walk that connects a vertex
$w\neq v'$ of $T(v)$ to a real vertex different from $v'$. Since
$T(v)$ is a tree there exists a function $f_v\in \mathcal{M} (T(v))$
such that $\Delta (f_v)=v-v'$ on $T(v)$. Let $f\in \mathcal{M} (G)$
such that $f|_{V(T(v))}=f_v$ and $f(v'')=f(v')$ for $v''\in
V(G)\setminus V(T(v))$. It follows that on $G$ we have $\Delta
(f)=v-v'$, hence $v\sim v'$. In the same way, using
$\overline{T(v)}$ we obtain that $\overline{v}\sim v'$. It follows
that $v+\overline{v}\sim 2v'$ on $G$.

Now assume $G(\mathbb{R})$ has $m>1$ connected components and assume
the theorem holds for M-graphs having less than $m$ connected
components. As in step 5 of the proof of Theorem \ref{theorem1} we
fix a component $G_1$ of $G(\mathbb{R})$ and let $G_2$ be as in that
proof. We make use of the results in the claim mentioned in step 5
of the proof of Theorem \ref{theorem1}. We obtain $G_2\subset G'_2$,
in particular $x\geq 2$, hence $G'_1, \cdots G'_x$ are graphs with a
real structure such that $G'_i(\mathbb{R})$ has less than $m$
connected components (say $G'_i(\mathbb{R})$ has $m_i$ connected
components and $\sum_{i=1}^x m_i=m$). Also $s(G'_i)=g(G'_i)+1$ for
$1\leq i\leq x$, hence $G'_1; \cdots; G'_x$ are M-graphs and we can
apply the induction hypothesis on them.

\begin{figure}[h]
\begin{center}
\includegraphics[height=6 cm]{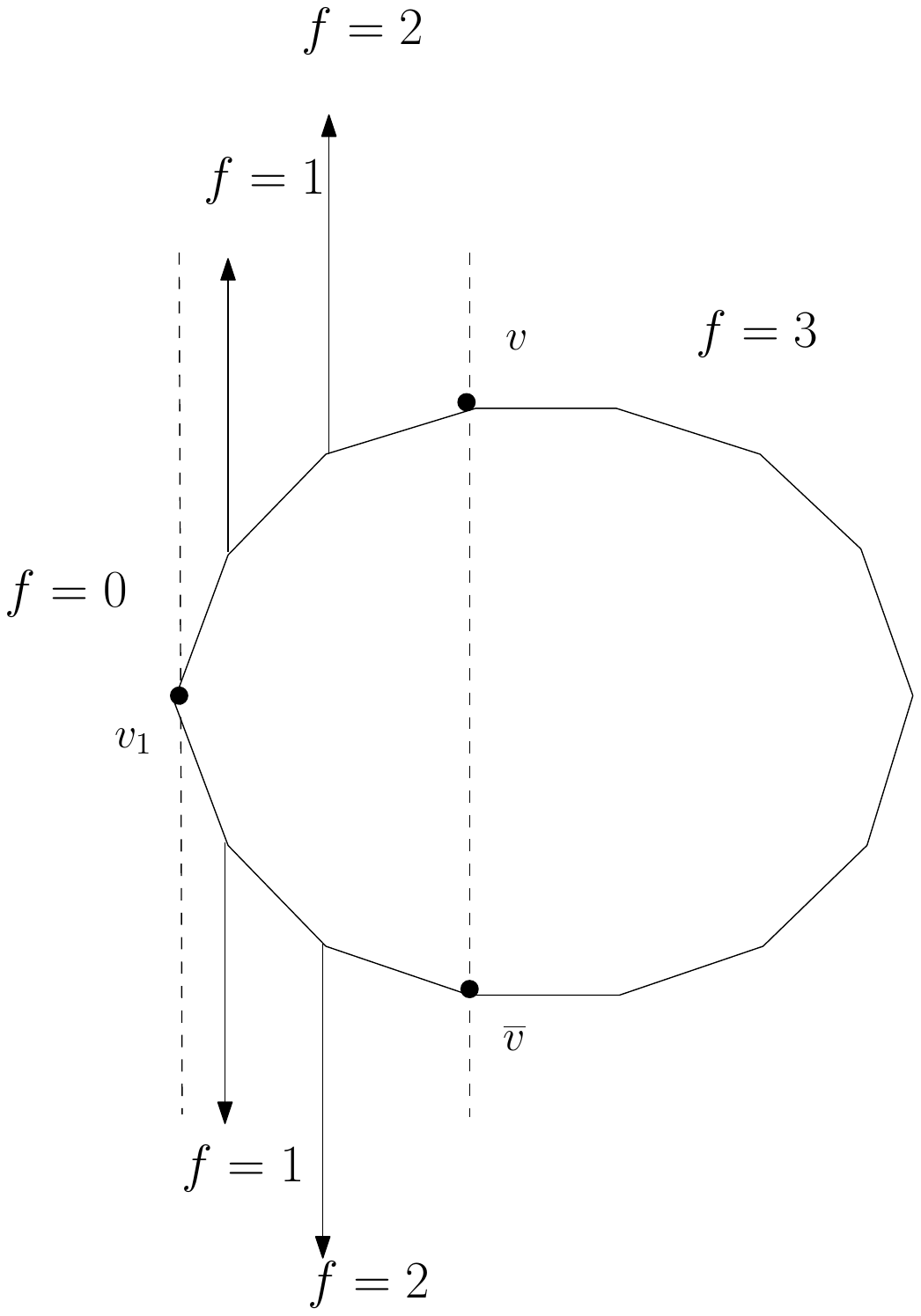}
\caption{$v+\overline{v}$ on $S$.}
\end{center}
\end{figure}

Let $v+\overline{v}$ be a non-real vertex pair on $S$ (the subgraph
used in step 5 of the proof of Theorem \ref{theorem1}). From the
above mentioned claim it follows that $S$ is a cycle. Let $\Gamma$
be the shortest walk on $S$ from $v_1$ to $v$; let it be
$v_1=w_0e_0w_1e_1 \cdots e_{a-2}w_{a-1}e_{a-1}w_a=v$. Let
$f(w_i)=f(\overline{w_i})=i$ for $0\leq i\leq a$ and $f(w)=a$ for
$w\in V(S)\setminus \{w_0; \cdots ; w_a; \overline{w_0}; \cdots ;
\overline{w_a}\}$. Let $f=0$ on $G'_1$ and $f=a$ on $G'_2$. For
$3\leq i\leq t$ let $f=f(v_i)$ on $G_i$ with $W_i=\{ v_i;
\overline{v_i}\}$ in case $i\leq x$ and $W_i=\{ v_i\}$ in case $i>x$
(see again the above mentioned claim). Then $v+\overline{v}+\Delta
(f)=2v_1$, hence $v+\overline{v}\sim 2v_1$. In a similar way one
finds $v+\overline{v}\sim 2v_2$, hence $2v_1\sim 2v_2$.

\begin{figure}[h]
\begin{center}
\includegraphics[height=4 cm]{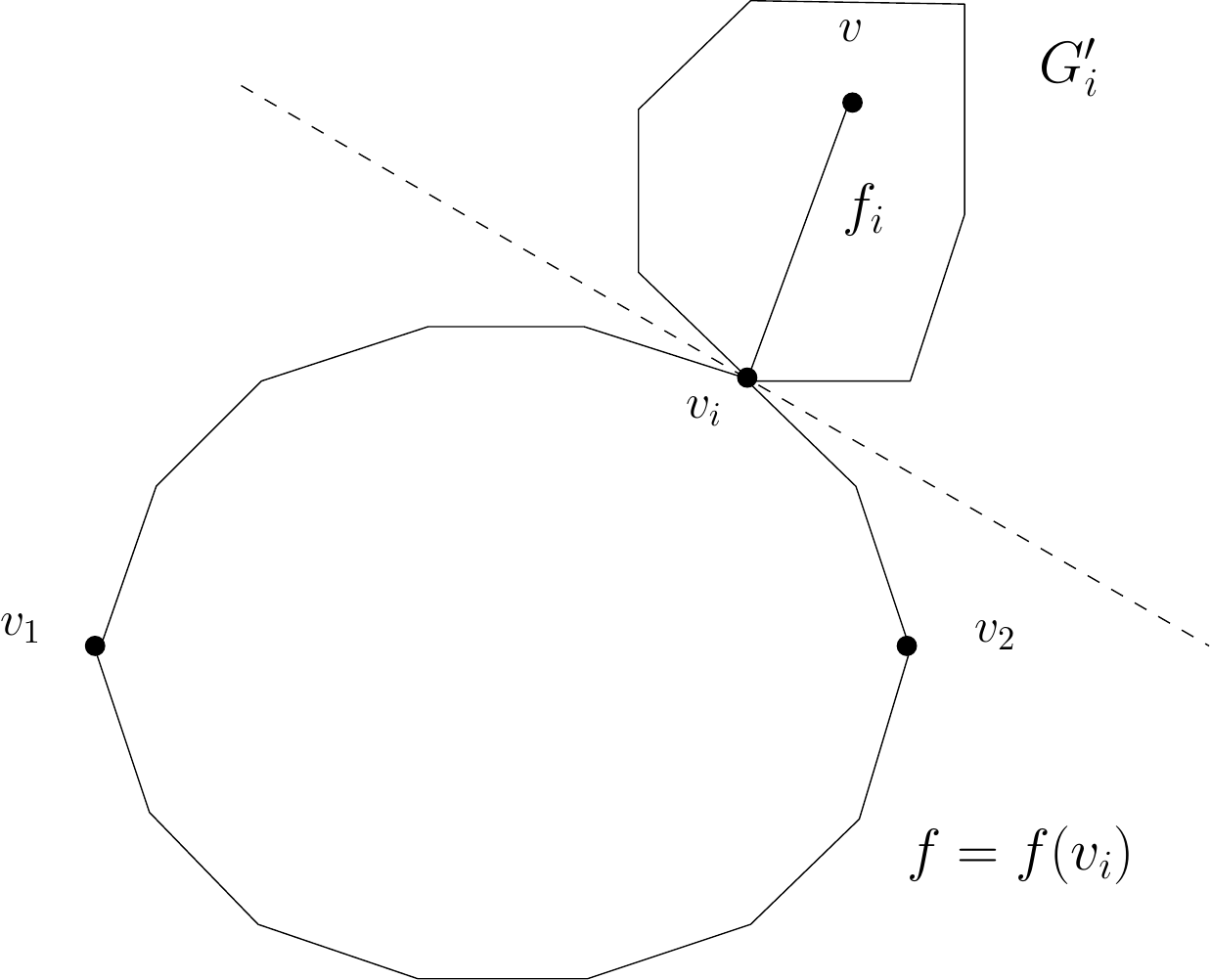}
\caption{$v$ on $G_i$ with $i>x$.}
\end{center}
\end{figure}

Next assume $v$ is a vertex of some $G'_i$ different from $v_i$ (as
above) for some $x+1\leq i\leq t$. As in the case with
$G(\mathbb{R})$ being connected we obtain a tree $T(v)$ contained in
$G'_i$ satisfying $T(v)\cap C=\{ v_i\}$ and no vertex $w$ of $T(v)$
different from $v_i$ is the end point of an edge not contained in
$T(v)$. Since $T(v)$ is a tree there exists $f_v\in
\mathcal{M}(T(v))$ with $\Delta (f_v)=v_i-v$ on $T(v)$. Take $f\in
\mathcal{M}(G)$ with $f|_{V(T(v))}=f_v$ and $f(v')=f(v_i)$ for
$v'\in V(G)\setminus V(T(v))$ then $\Delta (f)=v_i-v$ on $G$. Using
$\overline{f}$ we obtain $\overline{v}\sim \overline{v_i}$ hence
$v+\overline{v}\sim v_i+\overline{v_i}$. But $v_i+\overline{v_i}\sim
2v_1$ hence $v+\overline{v}\sim 2v_1$.

\begin{figure}[h]
\begin{center}
\includegraphics[height=4 cm]{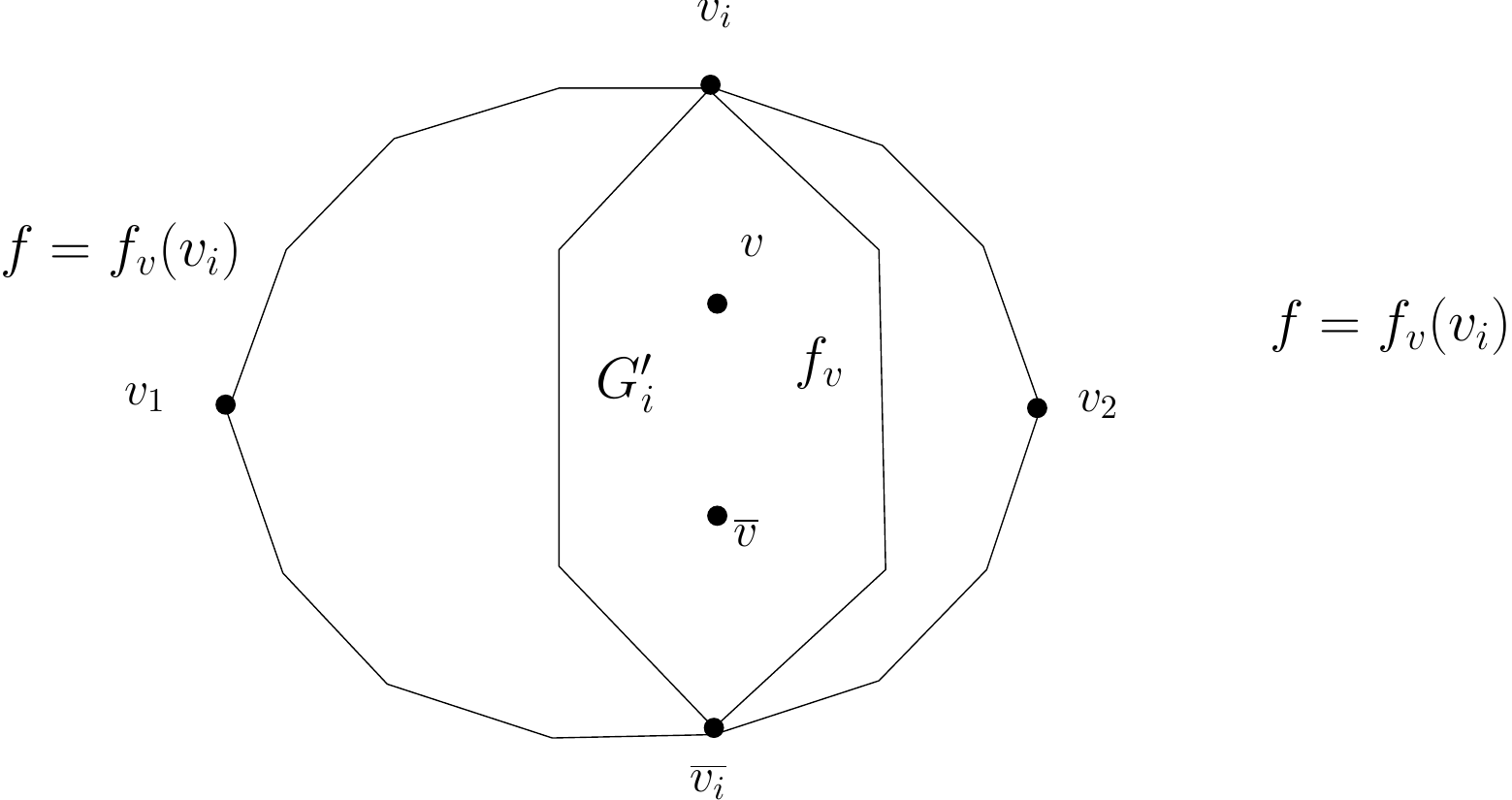}
\caption{$v+\overline{v}$ on $G_i$ with $i\leq x$.}
\end{center}
\end{figure}

Finally assume $v+\overline{v}$ is a non-real vertex pair in some
$G'_i$ with $1\leq i\leq x$. By induction there exists $f_i\in
\mathcal{M}(G'_i)$ and a real vertex $w$ of $G'_i$ such that $\Delta
(f_i)=2w-(v+\overline{v})$ on $G'_i$. From Lemma \ref{lemma1} we
know that $\overline{f_i}=f_i$, hence $f_i(v_i)=f_i(\overline{v_i})$
in case $i\geq 3$. Define $f\in \mathcal{M}(G)$ such that
$f|_{G'_i}=f_i$ and $f(v')=f_i(v_i)$ for $v'\in V(G)\setminus
V(G'_i)$. We obtain $\Delta (f)=2w-(v+\overline{v})$ on $G$ hence
$v+\overline{v}\sim 2w$.

\end{proof}

\begin{definition}\label{definitie10}
An M-graph $G$ is called a \emph{strong M-graph} if $G(\mathbb{R})$
has $g(G)+1$ connected components.
\end{definition}

In case $G$ is a strong M-graph it follows that each component $G_i$
of $G(\mathbb{R})$ is a tree. For strong M-graphs one has the
following strong result concerning linear systems.

\begin{proposition}\label{propositie3}
A strong M-graph $G$ has a real linear system $g^1_2$ (meaning it
has a real effective divisor $D$ of degree 2 with $rk(D)=1$).
\end{proposition}

\begin{proof}

We already know from Theorem \ref{theorem2} that for each non-real
vertex pair $v+\overline{v}$ there exists a real vertex $w$ with
$v+\overline{v}\sim 2w$. Therefore it is enough to prove the
following result: let $v'_1; v'_2$ and $v''_1; v''_2$ be two pairs
of real vertices, each pair belonging to one component of
$G(\mathbb{R})$, then $v'_1+v'_2\sim v''_1+v''_2$ on $G$.

First assume $G(\mathbb{R})$ has only one component. Since $G$ is a
strong M-graph it follows that $G$ is a tree, hence $v'_1+v'_2\sim
v''_1+v''_2$. So assume $G(\mathbb{R})$ has $m\geq 2$ components and
assume the proposition holds for strong M-graphs having less than
$m$ components. Of course we can assume $v'_1; v'_2$ belong to
$G_1$. If $v''_1; v''_2$ also belong to $G'_1$ then we can use the
induction hypothesis to $G'_1$ (in the proof of Theorem
\ref{theorem2} we found $G'_1$ is an M-graph and since each
connected component of $G'_1(\mathbb{R})$ has genus 0 it is a strong
M-graph). There exists $f_1\in \mathcal{M}(G'_1)$ with $\Delta
(f_1)=v'_1+v'_2-v''_1-v''_2$. Define $f\in \mathcal{M}(G)$ such that
$f|_{G'_1}=f_1$ and $f(v')=f(v_1)$ for $v'\in V(G)\setminus
V(G'_1)$. Then $\Delta (f)=v'_1+v'_2-v''_1-v''_2$ on $G$. From this
part of the proof it follows that $v'_1+v'_2\sim 2v_1$. But in the
proof of Theorem \ref{theorem2} we also found $2v_1\sim 2v_2$ and
$2v_1\sim v_i+\overline{v_i}$ for $3\leq i\leq x$. From this part of
the proof, if $v''_1; v''_2$ belong to $G'_2$ (resp. $G'_i$ with
$i\geq 3$) we also have $v''_1+v''_2\sim 2v_2$ (resp.
$v''_1+v''_2\sim v_i+\overline{v_i}$) hence $v''_1+v''_2\sim 2v_1$
in those cases too, hence $v'_1+v'_2\sim v''_1+v''_2$.

\end{proof}

\section{Extensions to metric graphs and tropical
curves}\label{section4}

In order to study the behavior of the specialization of linear
systems on curves to linear systems on graphs as described in
\cite{ref2} with respect to the occurrence of a real structure , we
need to extend the concept of graphs with a real structure to the
context of metric graphs. As a slightly further generalization we
generalize to the context of tropical curves. We also consider
generalizations of results on real linear systems on graphs with a
real structure to those contexts.

A \emph{weighted graph} $G$ is a graph $G$ having a weight function
$w : E(G)\rightarrow \mathbb{R}_{>0}$. For an edge $e \in E(G)$ the
number $w(e)$ is called the weight of $e$.

\begin{definition}\label{definition11}
A \emph{real structure} on a weighted graph $G$ is a real structure
on the underlying graph such that for each $e\in E(G)$ one has
$w(\overline{e})=w(e)$.
\end{definition}

A \emph{metric graph} $\Gamma$ is a compact, connected metric space
such that each point $p$ on $\Gamma$ has a neighborhood that is
isometric to a star-shaped set of some valence $n_p\geq 1$. A
\emph{star-shaped set of valence $n$} is a metric subspace $S$ of
$\mathbb{C}$ obtained as follows. There exists a real number $r>0$
such that

\begin{displaymath}
    S= \{ z\in \mathbb{C} : z=te^{2k\pi i/n} \text { with } 0\leq t<
    r \text { and } k\in \mathbb{Z} \}
\end{displaymath}

A discrete subset $V$ of a metric graph $\Gamma$ is called a set of
vertices of $\Gamma$ if it contains each point $p$ satisfying
$n_p\neq 2$. Once such a set is chosen, we denote it by $V(\Gamma)$
and we call it the set of vertices of $\Gamma$. The connected
components of $V\setminus V(\Gamma)$ are isometric to line segments,
they are called the \emph{edges} of $\Gamma$. We write $E(\Gamma)$
to denote the set of edges of $\Gamma$. We write $\psi (e)$ to
denote the set of the end points of the closure $\overline{e}$ (i.e.
$\psi (e)=\overline{e}\cap V(\Gamma)$). The length of an edge $e$ of
$\Gamma$ is denoted by $l(e)$. Since $\Gamma$ is compact those
lengths $l(e)$ are finite real numbers.

\begin{definition}\label{definitie12}
A \emph{real structure} on a metric graph $\Gamma$ is an isometry
$\iota :\Gamma \rightarrow \Gamma$ such that $\iota ^2$ is the
identity.
\end{definition}

Let $\Gamma$ be a metric graph with a real structure. We can assume
(and we do) that the set of vertices is invariant under $\iota$. As
usual, for $p\in \Gamma$ we write $\overline{p}$ instead of $\iota
(p)$ and we say $\overline{p}$ is the \emph{conjugated point}. In
case $p=\overline{p}$ we say $p$ is \emph{real point} on $\Gamma$.
For $e\in E(G)$ the image $\iota (e)$ is an edge $\overline{e}\in
E(G)$. We say $e$ is a \emph{real edge} if $e$ is pointwise fixed by
$\iota$. In case $e=\overline{e}$ ($e$ is fixed by $\iota$) either
$e$ is real or conjugation on $e$ is given by reflection on $e$ with
center the mid-point of $e$. A real point on $\Gamma$ is either a
real vertex of $\Gamma$, a point on a real edge of $\Gamma$ or a
mid-point of a non-real edge $e$ of $\Gamma$ fixed by $\iota$.

Let $G$ be a weighted graph. Associated to $G$ there is a metric
graph $\Gamma (G)$ with bijections $b_V : V(G) \rightarrow V(\Gamma
(G))$ and $b_E : E(G) \rightarrow E(\Gamma (G))$ such that for $e\in
E(G)$ one has $l(b_E(e))=w(e)$ and $\psi (b_E(e))=b_V(\psi (e))$.

\begin{definition/construction}\label{definitie13}
Let $G$ be a weighted graph with a real structure. On the associated
metric graph $\Gamma (G)$ we define the associated real structure as
follows. We use the same conjugation on $V(\Gamma (G))$ as on $V(G)$
(using $b_V$). Let $e\in E(G)$ be an edge of $G$. If $e$ is an
non-isolated real edge of $G$ then $\iota |_{b_E(e)}$ is the
identity and in case $e$ is an isolated real edge $\iota |_{b_E(e)}$
is the reflection on $b_E(e)$ with center the mid-point of $b_E(e)$.
In case $e\in E(G)$ is not real with $\psi (e)=\{ v_1; v_2 \}$ (of
course $v_1=v_2$ is possible) we use an isometry $b_E(e) \rightarrow
b_E(\overline{e})$ mapping $b_V(v_i)$ to $b_V(\overline{v_i})$. This
isometry is unique if $v_1\neq v_2$ but there are 2 choices if
$v_1=v_2$. In the case when the isometry is not unique the choices
are made in a compatible way in order to obtain a real structure.
\end{definition/construction}

\begin{remark}
In the specialization described in Section \ref{section5} the graph
$G$ has no loops, hence the problem of the possible choices in
Definition/Construction \ref{definitie13} does not occur.
\end{remark}

Let $\Gamma$ be a metric graph. Then there is an associated weighted
graph $G(\Gamma)$ obtained by identifying $V(G(\Gamma))$ with
$V(\Gamma)$ and $E(G(\Gamma))$ with $E(G)$ and using the same
function $\psi$ for $G(\Gamma)$ as for $\Gamma$. Moreover for $e\in
E(G(\Gamma))$ we put $w(e)=l(e)$. Again a real structure on $\Gamma$
induces a (now uniquely defined) real structure on $G(\Gamma)$ using
the action $\iota$ on $V(G)$ and $E(G)$. In case $e$ is a non-real
edge of $\Gamma$ fixed by $\iota$ the edge $e$ of $G(\Gamma)$ is an
isolated real edge.

Since a weighted graph with a real structure is a graph with a real
structure Theorem \ref{theorem1} holds in the context of weighted
graphs.

On a metric graph $\Gamma$ with a real structure we write $\Gamma
(\mathbb{R})$ to denote the set of real points on $\Gamma$. It has a
finite number of connected components $\Gamma_1; \cdots;
\Gamma_{s'}$ each one being a metric subgraph (some of them can be
points). Let $s(\Gamma)=\sum_{i=1}^{s'}(g(\Gamma_i)+1)$. A
\emph{path} $P$ in a metric graph $\Gamma$ is the image of a
continuous map $\gamma :[0; 1]\rightarrow \Gamma$ and we say $P$
connects $\gamma (0)$ to $\gamma (1)$. We define $a(\Gamma)=1$ if
there exists a non-real point $p$ on $\Gamma$ such that $p$ is
connected to $\overline{p}$ by a path $P$ not containing any real
point of $\Gamma$, otherwise we define $a(\Gamma)=0$. Clearly
$g(\Gamma)=g(G(\Gamma))$, $s(\Gamma)=s(G(\Gamma))$ and
$a(\Gamma)=a(G(\Gamma))$, hence we obtain the following result on
metric graphs with a real structure.

\begin{proposition}\label{proposition4}
If $\Gamma$ is a metric graph with a real structure then
$s(\Gamma)\leq g(\Gamma)+1$; $s(\Gamma) \equiv g(\Gamma)+1 \pmod{2}$
and $s(\Gamma)\leq g(\Gamma)-1$ if $a(\Gamma)=1$.
\end{proposition}

\begin{remark}
This proposition on metric graphs with a real structure can also be
obtained as a special case of so-called Smith Theory in topology.
The inequality $s(\Gamma)\leq g(\Gamma)+1$ follows from
\cite{refT}*{Theorem 4.1} while the equality $s(\Gamma) \equiv
g(\Gamma)+1 \pmod{2}$ follows from \cite{refT}*{Theorem 4.3}. From
this one obtains the similar statements on graphs in Theorem
\ref{theorem1}. Also the structure on the M-graphs needed to prove
Theorem \ref{theorem3} and Proposition \ref{propositie3} can be
deduced from it. The author likes to thank prof. V. Kharlamov for
mentioning this relation to Smith Theory.
\end{remark}

A \emph{tropical curve} $T$ is a connected metric space being the
union of a metric graph $\Gamma$ and a finite number of unbounded
edges $e$. Such an unbounded edge is isometric to $[0; \infty]$ and
satisfies $\overline{T\setminus e}\cap e$ is a unique point on
$\Gamma$ corresponding to 0 and denoted by $e(0)$. Then $\psi (e)$
are the points on $e$ corresponding to 0 and $\infty$. We call
$\Gamma$ the \emph{finite part} of $T$.

\begin{definition}\label{definition14}
A \emph{real structure} on a tropical curve $T$ is an isometry
$\iota :T \rightarrow T$ such that $\iota ^2$ is the identity.
\end{definition}

A real structure on a tropical curve $T$ induces a real structure on
the finite part $\Gamma$. It also induces an involution on its set
of unbounded edges. It $e$ is such an unbounded edge with
$e=\overline{e}$ then $e(0)=\overline{e(0)}$ and $e$ is pointwise
fixed by $\iota$. In that case $e$ is a real unbounded edge of $T$.
As before a point $p$ on $T$ such that $p=\overline{p}$ is called a
\emph{real point} and a real point is either a real vertex of $T$, a
point on a real edge of $T$ or a mid-point on a non-real edge of
$\Gamma$ fixed by $\iota$. The set $T(\mathbb{R})$ has finitely many
connected components $T_1; \cdots ; T_{s'}$ each one being a
tropical curve with $T_i \cap \Gamma\neq \emptyset$. As in the case
of a metric graph we define numbers $s(T)$ and $a(T)$ and clearly
$g(\Gamma )=g(T)$, $s(\Gamma )=s(T)$ and $a(\Gamma )=a(T)$. This
implies the following proposition.

\begin{proposition}\label{proposition5}
Let $T$ be a tropical curve with a real structure, then $s(T)\leq
g(T)+1$, $s(T) \equiv g(T)+1 \pmod{2}$ and $s(T)\leq g(T)-1$ if
$a(T)=1$.
\end{proposition}

Now we extend the results of Section \ref{section3} to the context
of metric graphs and tropical curves. First we recall the definition
of linear systems on metric graphs and tropical curves (see e.g.
\cite{ref7}).

Let $\Gamma$ be a metric graph or more general a tropical curve. A
divisor on $\Gamma$ is a finite $\mathbb{Z}$-linear combination $D$
of points on $\Gamma$. For a point $p$ on $\Gamma$ we denote $D(p)$
for the coefficient of $D$ at $p$. A \emph{rational function} $f$ on
$\Gamma$ is a continuous mapping $f:\Gamma \rightarrow \mathbb{R}$
such that for each edge $e$ of $\Gamma$ (including the unbounded
edges if $\Gamma$ is a tropical curve) identified isometrically with
an interval $I\subset [0; +\infty]$ there is a finite partition $e_1
\cup \cdots \cup e_n$ of $e$ in subintervals such that $f| _{e_i}$
is affine with an integer slope. We write $\mathcal{M}(\Gamma)$ to
denote the set of rational functions on $\Gamma$.

In case $p\in \Gamma$ and $f\in \mathcal{M}(\Gamma)$, $\Delta
(f)(p)$ is the sum of slopes of $f$ on $\Gamma$ in all directions
emanating from $p$. In case $p$ is a vertex of $\Gamma$ and $e$ is
an edge of $\Gamma$ with $p\in \psi (e)$ there is one such slope
associated to $e$ at $p$ denoted by $s_e(f,p)$. In case $p\in
\Gamma$ is not a vertex there is a unique edge $e$ of $\Gamma$ with
$p\in e$ and there are two such slopes associated to $e$ at $p$
denoted by $s'_e(f,p)$ and $s''_e(f,p)$. In case $f$ is affine at
$p$, $s'_e(f,p)+s''_e(f,p)=0$. In this way we define a principal
divisor $\Delta (f)=\sum_{p\in \Gamma}\Delta (f)(p)$ on $\Gamma$.
Now we can define linear equivalence of divisors on $\Gamma$, linear
systems on $\Gamma$ and we can define the rank of a divisor on
$\Gamma$ as we did in the case of graphs.

In case $\Gamma$ has a real structure, as in the case of graphs for
a divisor $D$ on $\Gamma$, there is a conjugated divisor
$\overline{D}$ and for a rational function $f$ on $\Gamma$ there is
a conjugated rational function $\overline{f}$ on $\Gamma$. A divisor
$D$ is called a real divisor if $D=\overline{D}$ and it is called
totally real if moreover the support of $D$ is contained in $\Gamma
(\mathbb{R})$.

A metric graph $\Gamma$ is called \emph{rational} if the length of
each edge of $\Gamma$ is a rational number. In such case
$\Gamma_{\mathbb{Q}}$ is the set of points $p$ of $\Gamma$ such that
there is an edge $e\in E(\Gamma)$ with $p\in e$ and the distance of
$p$ to the points of $\psi (e)$ is a rational number. We write $\Div
^{\mathbb{Q}}(\Gamma)$ to denote the divisors having support on
$\Gamma_{\mathbb{Q}}$. If moreover $\Gamma$ has a real structure
then $\Div ^{\mathbb{Q}}_{\mathbb{R}}(\Gamma)$ is the group of real
rational divisors on $\Gamma$.

Some results on linear systems proved for graphs with a real
structure also do hold in the context of metric graphs and tropical
curves.

\begin{proposition}\label{proposition6}
Let $\Gamma$ be a metric graph or a tropical curve with a real
structure and let $D$ be a real divisor on $\Gamma$ with $\dim | D |
\geq r\geq 1$. Then for each effective real divisor $E$ on $\Gamma$
with $\deg (E)=r$ there exists a real divisor $D'\in | D |$ with
$D'\geq E$.
\end{proposition}

\begin{proof}

The proof is completely similar to the proof of Proposition
\ref{propositie1}. There exists $f\in \mathcal{M}(\Gamma)$ and
$D''\in | D |$ with $D''=D+ \Delta (f)$ and $D'' \geq E$. Let
$g=\max \{ f, \overline{f} \}$ and let $D'=D+ \Delta (g)$. Since
$\overline{g} =g$ and $D$ is real also $D'$ is real. Take $p\in
\Gamma$ and let $e$ be an edge of $\Gamma$ with $p\in e$. If
$f(p)\geq \overline{f}(p)$ then $g(p)=f(p)$ hence $s_e(g,p)\geq
s_e(f,p)$ in case $p$ is a vertex and $s'_e(g,p)\geq s'_e(f,p)$,
$s''_e(g,p)\geq s''_e(f,p)$ in case $p$ is not a vertex. If follows
that $D'(p)\geq D''(p)$ in that case. Similarly, in case
$\overline{f} (p)\geq f(p)$ we obtain $D'(p)\geq D''(\overline{p})$.
This implies $D'$ is effective and $D'\geq E$.

\end{proof}

\begin{lemma}\label{lemma2}
Assume $D_1$ and $D_2$ are real divisors on a metric graph or a
tropical curve $\Gamma$ with a real structure and assume $f\in
\mathcal{M}(\Gamma)$ such that $\Delta (f)=D_2-D_1$ then
$\overline{f}=f$.
\end{lemma}

The proof of this lemma is exactly the same as the proof of Lemma
\ref{lemma1}.

\begin{theorem}\label{theorem4}
Let $\Gamma$ be a metric graph or a tropical curve with a real
structure and let $D_1, D_2$ be two linearly equivalent real
divisors on $\Gamma$. If $\Gamma '$ is a connected component of
$\Gamma (\mathbb{R})$ then $\deg (D_1 |_{\Gamma '}) \equiv \deg (D_2
|_{\Gamma '}) \pmod{2}$.
\end{theorem}

\begin{proof}
Let $\Delta (f)=D_1-D_2$, hence $\overline{f}=f$ by Lemma
\ref{lemma2}. First assume $\Gamma '$ is a mid-point $p$ of a
non-real edge $e$ invariant for $\iota$. Since $\overline{f}=f$ and
$\iota$ is reflection on $e$ with center $p$ it follows that
$s'_e(f,p)=s''_e(f,p)$ hence $\Delta (f)(p)$ is even. Hence $D_1(p)
\equiv D_2(p) \pmod{2}$.

Now assume $\Gamma '$ is a connected component of $\Gamma
(\mathbb{R})$ that is not an isolated point unless it is a real
vertex of $\Gamma$. Use $f'=f|_{\Gamma '}$, similar to $f'$ in the
proof of Theorem \ref{theorem2}. Take $p\in \Gamma '$. If $p$ is not
a vertex of $\Gamma$ then $\Delta (f')(p)=\Delta (f)(p)$. Since
$\deg (D_1|_{\Gamma '} +\Delta (f'))=\deg (D_1|_{\Gamma '})$ such
points cannot contradict the claim of the theorem. In case $p$ is a
vertex of $\Gamma$ one can copy the arguments given in Theorem
\ref{theorem2}.

\end{proof}

\begin{definition}\label{definition15}
A metric graph $\Gamma$ or a tropical curve $T$ with a real
structure is called an \emph{M-metric graph} or an \emph{M-tropical
curve} if $s(\Gamma)=g(\Gamma)+1$ or $s(T)=g(T)+1$. It is called a
\emph{strongly M-metric graph} or a \emph{strongly M-tropical curve}
if the number of connected components of $\Gamma (\mathbb{R})$ or
$T(\mathbb{R})$ is equal to $g(\Gamma)+1$ or $g(T)+1$.
\end{definition}

In this definition there is no need to exclude something similar to
isolated real edges as in the case for graphs in Definition
\ref{definitie8} because they became isolated real points similar to
the isolated real vertex as obtained in step 1 of the proof of
Theorem \ref{theorem1}. Since Proposition \ref{proposition4} is
obtained from Theorem \ref{theorem1} in a direct way using weighted
graphs, it follows that the structure of M-metric graphs and
M-tropical curves is very similar to the structure of M-graphs.
Since this structure is the basis in proving Theorem \ref{theorem3}
and Proposition \ref{propositie3} we obtain similar results in the
context of metric graphs and tropical curves. We leave the details
of the proofs to the reader. It consists of using suited functions
on $\Gamma$ having slope 1 at some part of $\Gamma$ and being
constant on the connected components of the complement.

\begin{theorem}\label{thoerem5}
Let $\Gamma$ be an M-metric graph or let $T$ be an M-tropical curve.
Each real effective divisor is linearly equivalent to a totally real
effective divisor.
\end{theorem}

\begin{proposition}\label{proposition7}
Let $\Gamma$ be a strong M-metric graph or let $T$ be a strong
M-tropical curve. It has a real linear system $g^1_2$.
\end{proposition}

\section{Specialization of real linear systems on real curves to
graphs}\label{section5}

In his paper \cite{ref2} Baker introduces a specialization of linear
systems from curves to graphs using a specialization from smooth
curves to so-called strongly stable curves. We are going to show
that this specialization has good behavior with respect to real
structures.

Let $T$ be a smooth geometrically irreducible curve defined over
$\mathbb{R}$ and let $o\in T(\mathbb{R})$. Let $\pi :\mathcal{C}
\rightarrow T$ be a proper flat morphism of relative dimension 1
such that the general fiber is a smooth geometrically irreducible
curve and the total space $\mathcal{C}$ is a smooth surface. The
starting point is a local description of this situation. Since $T$
contains a smooth $\mathbb{R}$-rational point the function field of
$T$ is a real field (this means -1 is not a square inside
$\mathbb{R}(T)$; see \cite{ref8}*{p. 282}). The local ring of $T$ at
$o$ is a discrete valuation ring $R_0$ having residue field
$\mathbb{R}$. It follows that the completion is a valuation ring $R$
with residue field $\mathbb{R}$ and the quotient field $Q=Q(R)$
again is a real field (see Lemma \ref{lemma4}). We use the base
extension $\Spec (R)\rightarrow T$ and we obtain an arithmetic
surface $\mathfrak{X}\rightarrow R$ (here we write $R$ instead of
$\Spec (R)$) having special fiber $X_0=\pi^{-1}(o)$ (a curve defined
over $\mathbb{R}$) and generic fiber $X$ defined over $Q$. We are
going to use some easy facts concerning $R\subset Q$.

\begin{lemma}\label{lemma3}
Let $Q$ be a real field containing $\mathbb{R}$ then $Q\otimes
_{\mathbb{R}}\mathbb{C}$ is a field.
\end{lemma}

\begin{proof}
Assume the $\mathbb{R}$-algebra $Q\otimes _{\mathbb{R}}\mathbb{C}$
is not a field and let $M$ be a maximal ideal of $Q\otimes
_{\mathbb{R}}\mathbb{C}$. Since both $Q$ and $\mathbb{C}$ are
naturally embedded into $Q\otimes _{\mathbb{R}}\mathbb{C}$ we obtain
embeddings $Q\subset Q\otimes _{\mathbb{R}}\mathbb{C}/M$ and
$\mathbb{C}\subset Q\otimes _{\mathbb{R}}\mathbb{C}/M$. But $\dim
_Q(Q\otimes _{\mathbb{R}}\mathbb{C})=2$ hence we obtain $Q\cong
Q\otimes _{\mathbb{R}}\mathbb{C}/M$. This would imply
$\mathbb{C}\subset Q$ contradicting the fact that $Q$ is a real
field.
\end{proof}

\begin{lemma}\label{lemma4}
Let $R$ be a complete discrete valuation ring containing
$\mathbb{R}$ such that its residue field is equal to $\mathbb{R}$.
Then the quotient field is a real field.
\end{lemma}

\begin{proof}
Assume $Q$ is not a real field. There exists a non-zero element
$f/g$ of $Q$ (with $f, g\in R$) such that its square is -1, hence
$f^2+g^2=0$ in $R$. Let $v$ be the valuation of $R$ and let $\pi \in
R$ be a uniformizing element (i.e. $v(\pi )=1$). Let $a=\min \{
v(f), v(g) \}$. There exists $f', g'\in \mathbb{R}$ such that
$v(f-f'\pi ^{a})$ and $v(g-g'\pi ^{a})$ are both at least equal to
$a+1$. At least one of those numbers $f'; g'$ is different from 0.
From $f^2 +g^2=0$ it follows that $f'^2+g'^2=0$. This is impossible,
hence a contradiction.
\end{proof}

\begin{lemma}\label{lemma5}
Let $R$ be a complete discrete valuation ring containing
$\mathbb{R}$ such that its quotient field $Q$ is a real field and
its residue field is a finite extension of $\mathbb{R}$. Then the
residue field is equal to $\mathbb{R}$.
\end{lemma}

\begin{proof}
Assume there exists a complete discrete valuation ring $R$
containing $\mathbb{R}$ such that its quotient field is a real field
$Q$ and the residue field is isomorphic to $\mathbb{C}$. There
exists $x\in R$ such that $v(x^2+1)>0$ (again $v$ is as in the proof
of Lemma \ref{lemma4}). Since $Q$ is a real field we know from Lemma
\ref{lemma3} that $Q_{\mathbb{C}}=Q\otimes _{\mathbb{R}}\mathbb{C}$
is a field. Let $R_{\mathbb{C}}$ be the complete discrete valuation
ring inside $Q_{\mathbb{C}}$ extending $R$. Since the residue field
of $R$ is equal to $\mathbb{C}$ it follows $R_{\mathbb{C}}$ is
completely ramified over $R$. Choose an uniformizing parameter $\pi$
in $R_{\mathbb{C}}$, then this implies each element of
$R_{\mathbb{C}}$ can be written in a unique way as a series
$\sum_{n=0}^\infty (a_n+b_nx)\pi ^n$ with $a_n, b_n\in \mathbb{R}$.
Such element belongs to $R$ if and only if $a_n=b_n=0$ for odd $n$.
Since $i\in R_{\mathbb{C}}$ we have such an expression
$i=\sum_{i=0}^{\infty}(a_n+b_nx)\pi ^n$. Also
$x^2=-1+\sum_{j=1}^{\infty}(\alpha _j+\beta _jx)\pi^{2j}$. An
expression $-1=\sum _{n=0}^{\infty}(x_n+y_nx)\pi ^n$ implies
$x_0=-1$; $y_0=0$ and $x_n=y_n=0$ for $n\geq 1$.

Now we use

\begin{equation}\label{equationX}
    i^2=-1=\sum_{n=0}^{\infty}(a_n+b_nx)^2\pi^{2n}+2\sum_{
    \substack{ i;j=0\\ i<j}}^{\infty}(a_i+b_ix)(a_j+b_jx)\pi^{i+j}.
\end{equation}

\noindent Writing Equation \ref{equationX} as $-1=\sum
_{n=0}^{\infty}(x_n+y_nx)\pi ^n$ we call $x_n+y_nx$ the coefficient
of $\pi^n$ in Equation \ref{equationX}.

The only contribution to $\pi^0$ in Equation \ref{equationX} comes
from the term

\[
    (a_0+b_0x)^2=a_0^2-b_0^2+2a_0b_0x+b_0^2\sum_{j=1}^{\infty}(\alpha_j+\beta_jx)\pi^{2j}.
\]

\noindent Hence this term gives no contribution to $\pi^l$ for some
odd $l$ and moreover $a_0^2-b_0^2=-1$ and $2a_0b_0=0$ implying
$a_0=0$ and $b_0=1$. More generally the terms $(a_n+b_nx)^2\pi^{2n}$
and $(a_i+b_ix)(a_j+b_jx)\pi^{i+j}$ with $i+j$ even do not give
contributions to $\pi^l$ for some odd $l$.

Assume $k\geq 0$ with $a_i=b_i=0$ for $i$ odd with $i\leq 2k-1$.
Then the only contribution to $\pi^{2k+1}$ comes from the term

\[
    2(a_0+b_0x)(a_{2k+1}+b_{2k+1}x)\pi^{2k+1} \equiv
    (-2b_{2k+1}+2a_{2k+1}x)\pi^{2k+1} \pmod{\pi^{2k+2}}.
\]

Since the coefficient of $\pi^{2l+1}$ in Equation \ref{equationX}
has to be equal to 0 one obtains $a_{2k+1}=b_{2k+1}=0$. This proves
$a_l=b_l=0$ for odd $l$, hence $i\in R$. This contradicts the fact
that $Q$ is a real field.

\end{proof}

\begin{remark}
This lemma does not hold in general without assuming completeness.
As an example, take the localization of $\mathbb{R}[X]$ at the
maximal ideal $<X^2+1>$. In this case there exist two extensions of
$R$ in $Q_{\mathbb{C}}$.
\end{remark}

To define the specialization we put some more conditions on the
family $\mathfrak{X}\rightarrow R$. From Lemma \ref{lemma3} we know
$Q\otimes _{\mathbb{R}}\mathbb{C}$ is a field. Let $R_{\mathbb{C}}$
be the complete valuation ring of $Q_{\mathbb{C}}$ extending $R$ and
let $\mathfrak{X}_{\mathbb{C}}\rightarrow R_{\mathbb{C}}$ be
obtained from the base change $R\subset R_{\mathbb{C}}$. The closed
fiber $X_{0,\mathbb{C}}$ is obtained from $X_0$ by making the base
extension $\mathbb{R}\subset \mathbb{C}$, hence there is a complex
conjugation on $X_{0,\mathbb{C}}$. We assume $X_{0,\mathbb{C}}$ is
\emph{strongly semistable}. This means we assume $X_{0,\mathbb{C}}$
is reduced, all singular points of $X_{0,\mathbb{C}}$ are nodes and
its irreducible components are smooth. Associated to
$X_{0,\mathbb{C}}$ there is a dual graph $G$ with a real structure
obtained as follows.

\begin{construction}\label{construction1}

\normalfont The set of vertices $V(G)$ corresponds to the set of
irreducible components of $X_{0,\mathbb{C}}$. For a vertex $v\in
V(G)$ we write $C_v$ for the corresponding component. The set of
edges $E(G)$ corresponds to the set of nodes of $X_{0,\mathbb{C}}$.
If $e$ is an edge of $E(G)$ then we also write $e$ to denote the
node of $X_{0,\mathbb{C}}$ and $\psi (e)=\{ v,w \}$ if and only if
$e\in C_v\cap C_w$. The complex conjugation on $X_{0,\mathbb{C}}$
induces an involution on both $V(G)$ and $E(G)$ and this defines a
real structure on $G$. As an example, an isolated real edge of $G$
corresponds to a node on $X_{0,\mathbb{C}}$ that is an isolated real
point.
\end{construction}

\begin{figure}[h]
\begin{center}
\includegraphics[height=4 cm]{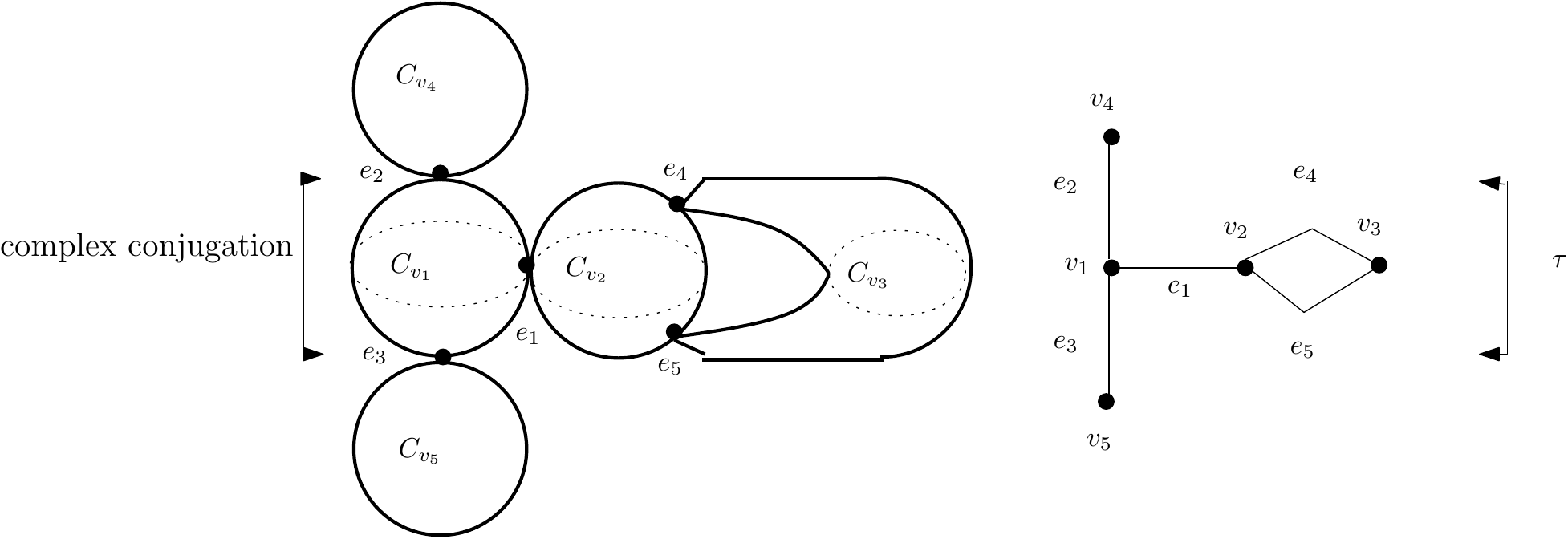}
\caption{Graph of a real curve}
\end{center}
\end{figure}

We write $\overline{Q}^r$ to denote the real closure of $Q$ and
$\overline{Q}$ to denote the closure of $Q$. Then
$\overline{Q}=\overline{Q}^r\otimes _{\mathbb{R}}\mathbb{C}$ (see
\cite{ref8}*{p. 274}). We write $\overline{X}^r$ (resp.
$\overline{X}$) to denote the curve obtained from $X$ by making the
base extension $Q\subset \overline{Q}^r$ (resp. $Q\subset
\overline{Q}$). On $\overline{Q}$ there is a natural involution
extending the complex conjugation on $\mathbb{C}$ and having
$\overline{Q}^r$ as its fixed field. We keep calling this complex
conjugation and it induces a so-called complex conjugation on
$\overline{X}$ too. The image of $P\in \overline{X}$ under this
complex conjugation is again called the conjugated point
$\overline{P}$. In case $P=\overline{P}$ we say $P$ is a real point
on $\overline{X}$. A divisor on $\overline{X}$ invariant under
conjugation is called a real divisor on $\overline{X}$.

Let $Q\subset K\subset \overline{Q}^r$ be an intermediate real field
with $[K:Q]=d$ finite. Let $R_K$ be the complete discrete valuation
ring extending $R$ having quotient field $K$. From Lemma
\ref{lemma5} we know the residue field of $R_K$ is equal to
$\mathbb{R}$. Consider the base change $\mathfrak{X} \times _R R_K
\rightarrow R_K$, let $R_{K,\mathbb{C}}$ be the complete discrete
valuation ring in $K \otimes _R \mathbb{C} =K_{\mathbb{C}}$ and
consider the base change $\left( \mathfrak{X} \times \times_R R_K
\right) _{\mathbb{C}} \rightarrow R_{K, \mathbb{C}}$. The special
fiber is defined over $\mathbb{R}$ and if $s$ is a node then it is
an $A_{d-1}$-singularity of $\left( \mathfrak{X} \times_R R_K
\right)_{\mathbb{C}}$. The union of those nodes is a 0-dimensional
subscheme of $\mathfrak{X} \times _R R_K$ (i.e. it is defined over
$\mathbb{R}$). The following lemma implies that the resolution of
singularities of all those $A_{d-1}$-singularities is defined over
$\mathbb{R}$.

\begin{lemma}\label{lemma6}
Let $X$ be a geometrically irreducible variety defined over
$\mathbb{R}$ and let $P\in X(\mathbb{C})\setminus X(\mathbb{R})$.
Let $\overline{X}$ be obtained from $X$ by using a base extension
$\mathbb{R}\subset \mathbb{C}$ and let $\widetilde{\overline{X}}$ be
the blowing-up of $\overline{X}$ at $P$ and $\overline{P}$, then
$\widetilde{\overline{X}}$ is defined over $\mathbb{R}$
\end{lemma}

\begin{proof}
We may assume $X$ is affine. Let $M_P$ (resp. $M_{\overline{P}}$) be
the maximal ideal corresponding to $P$ (resp. $\overline{P}$) and
let $I=M_P \cap M_{\overline{P}}$. Let $I_{\mathbb{R}}=\{ f\in I :
\overline{f}=f \}$ defining a closed subscheme of $X$. An element of
$I$ can be written as $f_1+if_2$ with $f_1, f_2$ defined over
$\mathbb{R}$. From $f_1(\overline{P})+if_2(\overline{P})=0$ it
follows $f_1(P)-if_2(P)=0$ hence $f_1(P)=f_2(P)=0$. This implies
$f_1, f_2 \in I_{\mathbb{R}}$. We show that the sheaf of ideals
induced by $I_{\mathbb{R}} \otimes _{\mathbb{R}} \mathbb{C}$ on
$\overline{X}$ is the sheaf of ideals induced by $I$. Clearly the
stalk of $I_{\mathbb{R}}\otimes _{\mathbb{R}}\mathbb{C}$ at $P$ is
included in $M_P$. In case $f\in I$, writing $f=f_1+if_2$ as before,
we have $f_1 \otimes 1+f_2 \otimes i\in I_{\mathbb{R}}\otimes
_{\mathbb{R}} \mathbb{C}$ maps to $f$. Assume $f\in M_P$ but
$f\notin M_{\overline{P}}$ hence $f(\overline{P})\neq 0$ and
therefore $\overline{f}(P)\neq 0$. Clearly $f\overline{f} \in
I_{\mathbb{R}}$. Since $\overline{f}$ is invertible at $P$ we can
considerer $\overline{f}^{-1}(f\overline{f} \otimes 1)$ mapping to
$f$. Hence the stalk at $P$ is equal to $M_P$. In a similar way the
stalk at $\overline{P}$ is equal to $M_{\overline{P}}$.

Let $Q$ be a point on $\overline{X}$ different from $P$ and
$\overline{P}$. If there exists $f\in M_P$ with
$f(Q)\overline{f}(Q)\neq 0$ then $f\overline{f} \otimes 1$ maps to
an invertible function at $Q$. Of course the same argument holds
using $\overline{P}$ instead of $P$. Assume for each $f\in M_P\cup
M_{\overline{P}}$ one has $f(Q)\overline{f}(Q)=0$. Since $Q\notin \{
P,\overline{P} \}$ there exists $f\in M_P$ and $f'\in
M_{\overline{P}}$ with $f(Q)f'(Q)\neq 0$. Since
$\overline{f}(Q)\overline{f'}(Q)=0$ and
$ff'+\overline{f}\overline{f'} \in I_{\mathbb{R}}$ it follows
$(ff'+\overline{f}\overline{f'}) \otimes 1$ maps to an invertible
function at $Q$. It follows that the stalk at $Q$ is the local ring
of $X$ at $Q$.

Now let $\widetilde{X} \rightarrow X$ be the blowing-up of $X$ along
$I_{\mathbb{R}}$. This is defined over $\mathbb{R}$. It follows from
e.g. \cite{ref9}*{B.6.9} that $\widetilde{\overline{X}}$ is a closed
subscheme of $\widetilde{X} \times _{\mathbb{R}} \mathbb{C}$. But
clearly $\widetilde{\overline{X}}$ is invariant under complex
conjugation in $\widetilde{X}\otimes _{\mathbb{R}} \mathbb{C}$ hence
$\widetilde{\overline{X}}$ is defined over $\mathbb{R}$.

\end{proof}

We write $\mathfrak{X}_K$ to denote the resulting smooth arithmetic
surface over $R_K$, it satisfies the same assumptions as
$\mathfrak{X} \rightarrow R$. In particular the special fiber of
$\mathfrak{X} _{K,\mathbb{C}} \rightarrow R_{K,\mathbb{C}}$ is a
strongly semistable curve $X_{0,K, \mathbb{C}}$ defined over
$\mathbb{R}$ and associated to it there is a graph $G_K$ with a real
structure.

\begin{construction}\label{construction2}

\normalfont Let $G_d$ be the graph obtained from $G$ by subdividing
each edge $e$ in $d$ parts. Given a real structure on $G$ we define
a real structure on $G_d$ as follows. There is a natural inclusion
$V(G)\subset V(G_d)$ and the conjugation on $V(G_d)$ restricts to
the known conjugation on $V(G)$. Let $e\in E(G)$ and $\psi (e)= \{
v; w \}$ (in general $v=w$ is possible). This edge is replaced by
edges $e_1; \cdots ; e_d$ with $\psi (e_i) = \{ v_{i-1}; v_i \}$
satisfying $v_0=v$ and $v_d=w$ (hence $v_i\notin V(G)$ for $1\leq
i\leq d-1$). First assume $\overline{e}\neq e$ then $\psi
(\overline{e})=\{ \overline{v}; \overline{w} \}$ and $\overline{e}$
is replaced by edges $\overline{e}_1; \cdots; \overline{e}_d$ with
$\psi (\overline{e}_i)= \{ \overline{v}_{i-1}; \overline{v}_i \}$
satisfying $\overline{v}_0=\overline{v}$ and
$\overline{v}_d=\overline{w}$. We define
$\overline{e_i}=\overline{e}_i$ and $\overline{v_i}=\overline{v}_i$
for $1\leq i\leq d$. In particular the edges $e_1; \cdots ; e_d$ and
the new vertices $v_1; \cdots; v_{d-1}$ are non-real. Assume $e$ is
a non-isolated real edge, hence also $v$ and $w$ are real vertices.
We define $\overline{v_i}=v_i$ and $\overline{e_i}=e_i$ for $1\leq
i\leq d$. Hence all edges $e_i$ and all new vertices $v_i$ are real.
Finally assume $e$ is an isolated real edge. Then $w=\overline{v}$.
We define $\overline{v_i}=v_{d-i}$ and $\overline{e_i}=e_{d+1-i}$.
Hence all edges $e_i$ and all new vertices $v_i$ are non-real except
for vertex $v_{d/2}$ in case $d$ is even and edge $e_{(d+1)/2}$ in
case $d$ is odd. In case $d$ is odd $e_{(d+1)/2}$ is an isolated
real edge.

\end{construction}

\begin{proposition}\label{proposition8}
The graph $G_K$ with its real structure is equal to the graph $G_d$
with its real structure.
\end{proposition}

\begin{proof}
The proper transforms of the components of $X_{0,\mathbb{C}}$ are
components of $X_{0,K,\mathbb{C}}$ and they give rise to a natural
inclusion $V(G)\subset V(G_K)$. Lemma \ref{lemma6} implies that the
conjugation on $V(G_K)$ induces the conjugation on $V(G)$. All
components corresponding to vertices from $V(G_K)\setminus V(G)$
correspond to new components obtained from resolving the
singularities of $\left( \mathfrak{X} \times _R R_K \right)
_{\mathbb{C}}$. Let $e\in E(G)$ with $\psi (e)= \{ v;w \}$. First
let $e$ be a non-real edge of $G$; it corresponds to a non-real
singular point (also denoted by $e$) of $\left( \mathfrak{X} \times
_R R_K \right) _{\mathbb{C}}$. Of course, the conjugated point
$\overline{e}$ is also a singular point of $\left( \mathfrak{X}
\times _R R_K \right) _{\mathbb{C}}$. On a small complex
neighborhood of $e$ on $\left( \mathfrak{X} \times _R R_K \right)
_{\mathbb{C}}$ we construct the resolution of the
$A_{d-1}$-singularity obtaining as exceptional divisor a chain of
$d-1$ copies of $\mathbb{P}^1_{\mathbb{C}}$ intersecting the proper
transforms of $C_v$ and $C_w$ transversally at one point at the ends
of the chain (those are the images on $C_v$ and $C_w$ of the
singularity $e$). In $V(G_K)$ this corresponds to subdividing $e$
into $d$ parts. From Lemma \ref{lemma6} it follows that using
complex conjugation on $\mathfrak{X}_{0,K,\mathbb{C}}$ we obtain the
chain of $\mathbb{P}^1_{\mathbb{C}}$'s obtained from resolving the
singularity of $\left( \mathfrak{X} \times _R R_K \right)
_{\mathbb{C}}$ at $\overline{e}$. On $V(G_K)$ we obtain the
subdividing of $\overline{e}$ into $d$ parts and the same
conjugation for this part of $V(G_K)$ as on the graph $G_d$.

Next assume $e$ is a non-isolated real edge of $G$. It corresponds
to a non-isolated real node of $X_{0,\mathbb{C}}$ and locally the
equation of $\mathfrak{X} \times _R R_K$ over $\mathbb{R}$ is given
by $x^2-y^2=t^d$. In this local description $x-y=0$ and $x+y=0$ are
local equations of $C_v$ and $C_w$. The blowing-up of $\mathbb{C}^3$
at $(0;0;0)$ and the strict transform $\mathfrak{X}'$ of the subset
of $\left( \mathfrak{X} \times _R R_K \right) _{\mathbb{C}}$ is
defined over $\mathbb{R}$. Let $X; Y; T$ be homogeneous coordinates
on the exceptional divisor $E$ then the intersection of $E$ with
$\mathfrak{X}'$ has equation $X^2-Y^2=T^2$ in case $d=2$ and
$X^2-Y^2=0$ in case $d\geq 3$. The intersection points with $C_v$
and $C_w$ are $(1:1:0)$ and $(1:-1:0)$. In case $d=2$ we obtain one
new component defined over $\mathbb{R}$. This corresponds to a new
real vertex $v_1$. In case $d\geq 2$ we obtain two new real
components corresponding to $X-Y=0$ and $X+Y=0$; those correspond to
the new vertices $v_1$ and $v_{d-1}$. The local equation of
$\mathfrak{X}'$ at $(0:0:1)$ is given by $x^2+y^2=t^{d-2}$. Hence if
$d=3$ the resolution is finished, if $d\geq 4$ we have to continue.
It follows that in $G_K$ the edge $e$ is divided into $d$ parts and
all new edges and vertices are real. This corresponds to the real
structure on $G_d$.

Finally assume $e$ is an isolated real edge. It corresponds to an
isolated real node of $X_{0,\mathbb{C}}$ and locally over
$\mathbb{R}$ the equation of $\left( \mathfrak{X} \times _R R_K
\right) _{\mathbb{C}}$ is given by $x^2+y^2=t^d$. Using a similar
blowing-up the intersection of $\mathfrak{X}'$ and $E$ has equation
$X^2+Y^2=T^2$ if $d=2$ and $X^2+Y^2=0$ if $d>2$. The intersection
points with the proper transform $C_v$ and $C_w=C_{\overline{v}}$
are $(1:i:0)$ and $(1:-i:0)$. Hence they are two conjugated non-real
nodes. This proves $\overline{e_1}=e_d$, hence $e_1$ (and $e_d$) are
non-real edges. In case $d=2$ there is a new component $v_1$ and it
is defined over $\mathbb{R}$ (it is the empty conic over
$\mathbb{R}$). In case $d\geq 3$ we obtain two conjugated components
with equations $X+iY=0$ and $X-iY=0$. This shows
$\overline{v_1}=v_{d-1}$. The local equation over $\mathbb{R}$ of
$\mathfrak{X}'$ at $(0:0:1)$ is given by $x^2+y^2=t^{d-2}$. In case
$d=3$ the singularity is resolved and the intersection of the new
components is an isolated real node corresponding to an isolated
real edge on $G_K$. In case $d\geq 4$ we have to continue. It
follows that in $G_K$ the edge $e$ is divided into $d$ parts and the
real structure on it corresponds to $G_d$.

\end{proof}

Now we recall the specialization map studied by Baker in
\cite{ref2}. For a divisor $\mathcal{D}$ on
$\mathfrak{X}_{\mathbb{C}}$ and a component $C$ of
$X_{0,\mathbb{C}}$ consider the intersection number defined by
$(C.\mathcal{D})= \deg (\mathcal{O}_{\mathfrak{X}}(\mathcal{D})|
_{C})$. This defines a homomorphism

\begin{displaymath}
\rho : \Div (\mathfrak{X}) \rightarrow \Div (G) : \mathcal{D}
\rightarrow \sum_{v\in V(G)} (C_v.\mathcal{D})v.
\end{displaymath}

\noindent Now let $D$ be a divisor on $X$ (hence a divisor on
$\overline{X}$ defined over $Q$). Clearly $D=D_1-D_2$ with $D_1$ and
$D_2$ both effective divisors on $X$, hence closed subschemes of
$X\subset \mathfrak{X}$. Take the closures of $D_1$ and $D_2$ on
$\mathfrak{X}$; their difference is denoted by $cl(D)\in \Div
(\mathfrak{X})$. The composition with $\rho$ gives rise to a
homomorphism (also denoted by $\rho$)

\begin{displaymath}
\rho : \Div (X)\rightarrow \Div (G) : D \rightarrow \rho (cl(D))
\end{displaymath}

\noindent Since $D$ is defined over $Q$ the divisor $cl(D)$ is
defined over $\mathbb{R}$ and this implies $\rho (D)$ is a real
divisor on $G$. Let $\Gamma$ be the metric graph associated with the
weighted graph obtained from $G$ giving weight 1 to each edge $e\in
E(G)$. The real structure on $G$ induces a real structure on
$\Gamma$ as explained in Definition/Construction \ref{definitie13}.
So we consider $\Gamma$ as a rational metric graph with a real
structure obtained in that way. From the previous map Baker obtains
a natural homomorphism $r: \Div (\overline{X}) \rightarrow \Div
^{\mathbb{Q}} (\Gamma)$ (see the proof of Theorem \ref{theorem6}).
It is a fundamental observation made by Baker that this map
preserves linear equivalence. Now we show that this natural
homomorphism also preserves the real structures on $X$ and $\Gamma$.
The preparations made in this section imply that we can follow
directly the arguments from \cite{ref2}.

\begin{theorem}\label{theorem6}
If $D$ is a divisor on $\overline{X}^r$ then $r(D)$ is a real
divisor on $\Gamma$. Hence $r$ induces a natural homomorphism
$r:\Div (\overline{X}^r) \rightarrow \Div
^{\mathbb{Q}}_{\mathbb{R}}(\Gamma)$.
\end{theorem}

\begin{proof}

Let $D$ be any divisor on $\overline{X}^r$ (i.e. a real divisor on
$\overline{X}$). Let $Q\subset K \subset \overline{K}^r$ be a finite
extension such that each real point $P$ and each non-real pair of
points $P+\overline{P}$ occurring in $\Supp (D)$ is defined over
$K$. It follows that each point occurring in $\Supp (D)$ as a
divisor on $\overline{X}$ is defined over $K_{\mathbb{C}}$. Write
$d=[K:Q]$ and consider $\mathfrak{X}_K \rightarrow R_K$. The special
fiber $X_{0,K,\mathbb{C}}$ of $\mathfrak{X}_{K,\mathbb{C}}
\rightarrow R_{K, \mathbb{C}}$ is defined over $\mathbb{R}$ (this
follows from Lemmas \ref{lemma5} and \ref{lemma6}) and from
Proposition \ref{proposition8} it follows that its dual graph with
real structure is equal to $G_d$ with the real structure obtained
from $G$ as described in Construction \ref{construction2}.

We consider $G_d$ as a weighted graph giving weight $1/d$ to each
edge $e\in E(G_d)$. Then the associated metric graph of $G_d$ is
equal to $\Gamma$ and this identification is compatible with the
real structures on $G_d$ and $\Gamma$. On
$\mathfrak{X}_{K,\mathbb{C}}$ each point $P$ of $D$ corresponds to a
closed point on the generic fiber, hence its closure is a section
$cl(P)$ of $\mathfrak{X}_{K,\mathbb{C}} \rightarrow
R_{K,\mathbb{C}}$. It intersects exactly one component of
$X_{0,K,\mathbb{C}}$ corresponding to a vertex $v(P)$ of $G_d$. This
vertex corresponds to rational point $r_K(P)$ on $\Gamma$. Then
$\rho (D)$ on $\Div (G_K)=\Div (G_d)$ is the sum of those vertices
$v(P)$ using $P$ on $D$, in particular it is a real divisor on
$G_d$. On $\Gamma$ we obtain $r_K(D)\in \Div
^{\mathbb{Q}}_{\mathbb{R}}(\Gamma)$. It is observed by Baker that
this divisor on $\Gamma$ does not depend on the chosen field $K$,
hence it is the divisor $r(D)$ on $\Gamma$.

\end{proof}

As a final remark we show that each graph $G$ with a real structure
can be obtained from a degeneration $\mathfrak{X}\rightarrow R$. We
are going to show that there exists a strongly semistable curve
$X_0$ defined over $\mathbb{R}$ having its dual graph equal to a
given graph with a real structure. As a matter of fact we are going
to obtain such curve that is totally degenerated (meaning each
component of $X_{0,\mathbb{C}}$ is a $\mathbb{P}^1_{\mathbb{C}}$).
Then the claim follows from a deep theorem (see e.g.
\cite{ref2}*{Appendix B}).

\begin{proposition}\label{proposition9}
Let $G$ be a graph of genus $g$ without loops having a real
structure. There exists a totally degenerated curve $X$ of genus $g$
defined over $\mathbb{R}$ such that the associated graph with real
structure is equal to $G$.
\end{proposition}

\begin{proof}

A real vertex corresponds to a copy of the real projective line
$\mathbb{P}^1_{\mathbb{R}}$ defined over $\mathbb{R}$.

A pair of non-real vertices $v+\overline{v}$ should correspond to a
curve $Y$ defined over $\mathbb{R}$ such that
$Y(\mathbb{R})=\emptyset$ and $Y(\mathbb{C})$ is the disjoint union
of two copies of the projective line defined over $\mathbb{C}$ and
interchanged by complex conjugation. Such a curve is obtained as
follows. Let $L, \overline{L}$ be two non-real conjugated lines in
$\mathbb{P}^2_{\mathbb{R}}(\mathbb{C})$ and $S=L\cap \overline{L}$,
hence $S\in \mathbb{P}^2_{\mathbb{R}}(\mathbb{R})$. The blowing-up
of $\mathbb{P}^2_{\mathbb{R}}$ at $S$ is a surface defined over
$\mathbb{R}$. The union of the strict transforms of $L$ and
$\overline{L}$ is defined over $\mathbb{R}$. This is $Y$.

Assume $v+\overline{v}$ is a pair of non-real vertices and assume
there exist exactly $m_1$ isolated real edges $e$ and $2m_2$
non-real edges $e$ with $\psi (e)=\{ v, \overline{v} \}$. Associated
to it we define a scheme $Y_{m_1, m_2}$ defined over $\mathbb{R}$
such that $Y(\mathbb{C})$ has two irreducible components $Y_1$ and
$Y_2$ both isomorphic to $\mathbb{P}^1_{\mathbb{C}}$, interchanged
by complex conjugation and such that $Y_1\cap Y_2$ consists of nodes
of $Y(\mathbb{C})$, $m_1$ of them being isolated real points and
$m_2$ pairs of conjugated non-real points. In order to obtain this
curve we take $Y$ as before, we choose $m_1$ different pairs of
conjugated points $P_1; \overline{P_1}; \cdots ; P_{m_1};
\overline{P_{m_1}}$ (with $P_i\in Y_1$ and $\overline{P_i}\in Y_2$)
and $m_2$ different pairs of non-conjugated points $Q_{1,1};
Q_{2,1}; \cdots ; Q_{1,m_2}; Q_{2,m_2}$ with $Q_{i,j}\in Y_i$ and
$Q_{2,j}\neq \overline{Q_{1,j'}}$ for all $j, j'$. Choose an
embedding of $Y$ is some $\mathbb{P}^N$ defined over $\mathbb{R}$
such that all chosen points together with their complex conjugated
points span a linear subspace $\Lambda$ of dimension $2m_1+4m_2-1$
such that $\Lambda$ intersects $Y$ with multiplicity 1 at those
points and contains no other point of $Y$. Let $R_j$ be a general
real point on the line connecting $P_j$ and $\overline{P_j}$ and let
$R'_j$ be a general point on the line connecting $Q_{1,j}$ and
$Q_{2,j}$. Let $\Lambda '$ be the linear span of the points $R_j$
and the points $R'_j$ and $\overline{R'_j}$. The projection with
center $\Lambda '$ defines the scheme $Y_{m_1,m_2}$ defined over
$\mathbb{R}$.

Now we can construct the curves using suited identifications of
points on two schemes $Y_1$ and $Y_2$ defined over $\mathbb{R}$ as
follows. Let $P_{1,1}; \cdots ; P_{1,m_1}$ and $P_{2,1}; \cdots ;
P_{2,m_1}$ be different real points on $Y_1$ and $Y_2$. Let
$Q_{1,1}; \cdots ; Q_{1,m_2}$ and $Q_{2,1}; \cdots : Q_{2,m_2}$ be
different non-real points on $Y_1$ and $Y_2$ with
$\overline{Q_{i,j}}\neq Q_{i,j}$ (the first index refers to the
curve). Using suited embeddings of $Y_1$ and $Y_2$ defined over
$\mathbb{R}$ and using a Segre embedding of the product one obtains
in a similar way as before using a suited projection defined over
$\mathbb{R}$ a new curve $Y$ defined over $\mathbb{R}$ obtained from
$Y_1$ and $Y_2$ by identifying $P_{1,j}$ to $P_{2,j}$; $Q_{1,j}$ to
$Q_{2,j}$ and $\overline{Q_{1,j}}$ to $\overline{Q_{2,j}}$. In this
way one obtains $m_1$ new non-isolated real nodes and $m_2$ new
pairs of complex conjugated non-real nodes. They do correspond to
$m_1$ real edges and $m_2$ pairs of conjugated non-real edges.
Iterating this construction one obtains as dual graph the given
graph with its real structure.

\end{proof}


\begin{bibsection}
\begin{biblist}

\bib{ref1}{article}{
    author={Baker, M.},
    author={Norine, S.},
    title={Riemann-Roch and Abel-Jacobi theory on a finite graph},
    journal={Adv. Math.},
    volume={215},
    year={2007},
    pages={766-788},
}

\bib{ref2}{article}{
    author={Baker, M.},
    title={Specialization of linear systems from curves to graphs. With an appendix by Brian Conrad},
    journal={Algebra Number Theory},
    volume={2},
    year={2008},
    pages={613-788},
}

\bib{ref6}{article}{
    author={Baker, M.},
    author={Norine, S.},
    title={Harmonic morphisms and hyperelliptic graphs},
    journal={Inst. Math. Res. Not. IMRN},
    volume={15},
    year={2009},
    pages={2914-2955},
}


\bib{ref10}{book}{
    author={Bondy, J.A.},
    author={Murty, U.S.R.},
    title={Graph Theory},
    series={Graduate Texts in Mathematics},
    volume={244},
    year={2008},
    publisher={Springer-Verlag},
}

\bib{refT}{book}{
    author={Brendon, G.E.},
    title={Introduction to compact transformation groups},
    series={Pure and Applied Mathematics},
    volume={46},
    year={1972},
    publisher={Academic Press, Inc.},
}


\bib{ref5}{article}{
    author={Ciliberto, C.},
    author={Harris, J.},
    title={Real abelian varieties and real algebraic curves},
    booktitle={Lectures in real geometry (Madrid, 1994)},
    series={de Gruyter Exp. Math.},
    volume={23},
    year={1996},
    pages={167-256},
}

\bib{ref9}{book}{
    author={Fulton, W.},
    title={Intersection theory},
    series={Ergebnisse der Mathematik und ihrer Grenzgebiete},
    volume={2},
    publisher={Springer-Verlag},
    year={1984},
}

\bib{ref7}{article}{
    author={Gathmann, A.},
    author={Kerber, M.},
    title={A Riemann-Roch theorem in tropical geometry},
    journal={Math. Z.},
    volume={259},
    year={2008},
    pages={217-230},
}

\bib{ref3}{article}{
    author={Gross, B.H.},
    author={Harris, J.},
    title={Real algebraic curves},
    journal={Ann. Scient. Ec. Norm. Sup.},
    volume={14},
    year={1981},
    pages={157-182},
}

\bib{ref4}{article}{
    author={Huisman, J.},
    title={On the geometry of algebraic curves having many real
    components},
    journal={Revista Matem\'atica Completense},
    volume={14},
    year={2001},
    pages={83-92},
}

\bib{ref8}{book}{
    author={Lang, S.},
    title={Algebra},
    publisher={Allison-Wesley Publ. Co.},
    year={1965},
}




\end{biblist}
\end{bibsection}
\end{document}